\batchmode
\makeatletter
\def\input@path{{/Users/chloekiss/Dropbox/Matrix-bipartite-cluster/Tech_note/}}
\makeatother
\documentclass[onecolumn,12pt,draftclsnofoot]{IEEEtran}
\usepackage[T1]{fontenc}
\usepackage[latin9]{inputenc}
\usepackage{color}
\usepackage{mathtools}
\usepackage{amsmath}
\usepackage{amsthm}
\usepackage{amssymb}
\usepackage{graphicx}
\PassOptionsToPackage{normalem}{ulem}
\usepackage{ulem}
\usepackage[pdftex,unicode=true,
 bookmarks=true,bookmarksnumbered=true,bookmarksopen=true,bookmarksopenlevel=1,
 breaklinks=false,pdfborder={0 0 0},pdfborderstyle={},backref=false,colorlinks=false]
 {hyperref}
\hypersetup{pdftitle={Your Title},
 pdfauthor={Your Name},
 pdfpagelayout=OneColumn, pdfnewwindow=true, pdfstartview=XYZ, plainpages=false}

\makeatletter
\theoremstyle{plain}
\newtheorem{thm}{\protect\theoremname}[section]
\theoremstyle{definition}
\newtheorem{defn}[thm]{\protect\definitionname}
\theoremstyle{remark}
\newtheorem{rem}[thm]{\protect\remarkname}
\theoremstyle{plain}
\newtheorem{lem}[thm]{\protect\lemmaname}
\theoremstyle{plain}
\newtheorem{prop}[thm]{\protect\propositionname}
\theoremstyle{plain}
\newtheorem{cor}[thm]{\protect\corollaryname}
\theoremstyle{definition}
\newtheorem{example}[thm]{\protect\examplename}
\theoremstyle{remark}
\newtheorem{claim}[thm]{\protect\claimname}

\usepackage[caption=false,font=footnotesize]{subfig}
\usepackage{algorithm,algpseudocode}

\usepackage{amstext}
\usepackage{pdfsync}
\usepackage{amsthm}
\usepackage{graphicx}
\usepackage{epstopdf}

\usepackage{tikz}
\usetikzlibrary{arrows,shadows,petri}
\tikzset{
  every overlay node/.style={
    draw=white,anchor=north west,
  },
}
\usetikzlibrary{positioning, decorations.pathmorphing}

\usepackage{hyperref}
\@ifundefined{rangeHsb}{\usepackage{xcolor}}{}

\usepackage{lipsum,multicol,cuted,mathtools}

\usepackage{scrlayer-scrpage}
\clearpairofpagestyles
\makeatother

\providecommand{\claimname}{Claim}
\providecommand{\corollaryname}{Corollary}
\providecommand{\definitionname}{Definition}
\providecommand{\examplename}{Example}
\providecommand{\lemmaname}{Lemma}
\providecommand{\propositionname}{Proposition}
\providecommand{\remarkname}{Remark}
\providecommand{\theoremname}{Theorem}

\begin{document}
\title{Sufficient Condition on Bipartite Consensus of Weakly Connected Matrix-weighted
Networks}
\author{Chongzhi~Wang, Haibin~Shao,~\IEEEmembership{Member,~IEEE,} Ying~Tan,~\IEEEmembership{Fellow,~IEEE,}
Dewei~Li\thanks{C.W., H.S., and D.L. are with the Department of Automation, Shanghai
Jiao Tong University, Shanghai, China, e-mail: \protect\href{mailto:\%7Bczhwon,shore,dwli\%7D@sjtu.edu.cn}{\{czhwon,shore,dwli\}@sjtu.edu.cn}.
C.W. and Y.T. are with the Department of Mechanical Engineering, University
of Melbourne, Melbourne, Australia, email: \protect\href{mailto:\%7Bchongzhi.wang.1,\%20yingt\%7D@unimelb.edu.au}{\{chongzhi.wang.1, yingt\}@unimelb.edu.au}.}}
\IEEEpubid{}
\maketitle
\begin{abstract}
Recent advances in bipartite consensus on matrix-weighted networks,
where agents are divided into two disjoint sets with those in the
same set agreeing on a certain value and those in different sets converging
to opposite values, have highlighted its potential applications across
various fields. Traditional approaches often depend on the existence
of a positive-negative spanning tree in matrix-weighted networks to
achieve bipartite consensus, which greatly restricts the use of these
approaches in engineering applications. This study relaxes that assumption
by allowing weak connectivity within the network, where paths can
be weighted by semidefinite matrices. By analyzing the algebraic constraints
imposed by positive-negative trees and semidefinite paths, we derive
new sufficient conditions for achieving bipartite consensus. Our findings
are validated by numerical simulations.
\end{abstract}


\IEEEpeerreviewmaketitle{}

\section{Introduction}

\IEEEPARstart{A}{chieving} collective behaviors among a group of
identical agents has been a prominent research topic across various
domains, including social networks \cite{Castellano2009,Proskurnikov2016},
coupled oscillator systems in power grids \cite{Doerfler2013}, and
the control of groups of robots such as unmanned ground vehicles \cite{Kamel2020},
unmanned aerial vehicles \cite{Tahir2019}, and unmanned underwater
vehicles \cite{Liu2023}. In addition to the dynamics of each agent,
communication protocols play a crucial role in facilitating these
collective group behaviors \cite{olfati2004consensus,altafini2012consensus,jadbabaie2003coordination}.
Early scalar-weighted protocols commonly used a scalar quantity for
each pair of communicating agents. While this approach ensures consensus
among high-dimensional agents with minimal requirements on network
topology, it does not fully utilize the extra communication capacity
of the agents. In contrast, recent research extends this concept by
using a matrix\textemdash typically non-diagonal\textemdash for each
pair of agents \cite{tuna2016synchronization,TRINH2018415,Pan2019},
which reduces to the scalar-weighted case when each matrix is a scalar
multiple of the identity matrix. A network employing this matrix-based
protocol is called the matrix-weighted network. It is shown to capture
more complex, inter-dimensional interactions and lead to more intricate
group behaviors \cite{tuna2016synchronization,Ramirez2009,zhao2016localizability,Lee2016}.

Among these group behaviors, bipartite consensus represents a specific
distribution of the agents. When reaching bipartite consensus, agents
are divided into two groups based on their states: one group consists
of agents with identical steady-state behaviors, while the other group
comprises agents with opposing steady-state behaviors. Bipartite consensus
finds applications in various fields, including social networks \cite{Proskurnikov2016},
coupled oscillator systems \cite{Hong2011}, and multi-robot systems
\cite{Hu2015,Zong2019}. This concept is also crucial for understanding
the more general case of cluster consensus, as bipartite consensus
represents a special instance of cluster consensus involving two clusters.

Research has focused on the conditions under which a matrix-weighted
network can achieve bipartite consensus. Given that the matrix-weighted
network can be modeled as a linear time-invariant system characterized
by a matrix-valued graph Laplacian $L$, the asymptotic behavior of
bipartite consensus is influenced by the structure of the null space
of the system matrix, ${\bf null}(L)$. In \cite{su2019bipartite},
a necessary and sufficient condition for achieving bipartite consensus
was provided, based on the algebraic properties of $L$. In contrast,
\cite{wang2022characterizing} proposed a necessary condition from
the perspective of network structure. This condition involves the
concept of a nontrivial balancing set, which, among other things,
divides the agents into two groups that are negatively interacting,
i.e., having negative connections in between. The study \cite{wang2022characterizing}
further shows that for a network to reach bipartite consensus, the
communication graph must contain a unique nontrivial balancing set
among all possible divisions. Although these conditions are theoretically
elegant, they lack practical design guidelines, such as sufficient
conditions that engineering practitioners can use to design matrix
weights and network structures to achieve bipartite consensus without
being overly conservative. In \cite{TRINH2018415} and \cite{Pan2019},
sufficient conditions are presented that require the communication
graph to have either a positive or a positive-negative spanning tree.
We refer to matrix-weighted networks with such a spanning tree as
a strongly connected network.

These requirements however, are not always met in practical applications.
For example, \cite{zhao2016localizability} applied the matrix-weighted
protocol in a sensor network where the connections are weighted solely
by semidefinite matrices. This network, lacking a positive-negative
spanning tree, is referred to as a weakly-connected matrix-weighted
network. Similar cases appear in coupled oscillator systems \cite{tuna2016synchronization,tuna2019synchronization}
and opinion dynamics \cite{Ahn2020}. This motivates this work, which
aims to provide sufficient conditions to ensure bipartite consensus
on weakly-connected matrix-weighted networks.

The weakly-connected matrix-weighted network poses a unique challenge
on the analysis of the Laplacian null space. The block Laplacian matrix
defined on these networks exhibits complex spectral properties and
null space behavior that are not characterized by traditional graph
metrics. To address this, we introduce new connectivity concepts to
partition the graph into components \textendash{} the strongly connected
subgraphs spanned by positive-negative trees and the semidefinite
paths between them \textendash{} enabling a more targeted null space
analysis. By exploring various necessary or sufficient conditions,
we develop a theorem that establishes the relationship between the
Laplacian null space and the network's structural and weight configurations,
which further ensures bipartite consensus. This work provides a theoretical
framework that not only addresses the bipartite consensus phenomenon
but also offers broader insights for designing and analyzing matrix-weighted
networks.

The paper is organized as follows. Section 2 introduces notations
and preliminary results of matrix-weighted networks. Section 3 presents
the definition and properties of the continent, followed by the main
theorem on achieving bipartite consensus. Numerical examples are demonstrated
in Section 4 to validate our theorem. Section 5 gives concluding remarks.

\section{Preliminaries and Problem Formulation}

This section covers the necessary preliminaries, including notations,
graph theory, problem formulation, and the relevant propositions.

\subsection{Notations}

Let $\mathbb{R}$, $\mathbb{N}$ and $\mathbb{Z}_{+}$ be the set
of real numbers, natural numbers and positive integers, respectively.
For $n\in\mathbb{Z}_{+}$, we denote $\underline{n}=\left\{ 1,2\text{,}\cdots,n\right\} $.
For $x\in\mathbb{R}^{n}$, the symbol $||\cdot||$ denotes its Euclidean
norm. For a set $M$, $\text{card}(M)$ denotes its cardinality. For
a symmetric matrix $Q\in\mathbb{R}^{n\times n}$, we express its positive
(negative) definiteness with $Q\succ0$ ($Q\prec0$), and its positive
(negative) semidefiniteness with $Q\succeq0$ ($Q\preceq0$). In addition,
a matrix-valued sign function ${\bf sgn}(\cdot)$ is defined such
that 
\[
\begin{cases}
{\bf sgn}(Q)=+1, & Q\succeq0,Q\neq0\text{ or }Q\succ0\\
{\bf sgn}(Q)=-1, & Q\preceq0,Q\neq0\text{ or }Q\prec0\\
{\bf sgn}(Q)=0, & Q=0.
\end{cases}
\]
The notation $|\cdot|$ is adopted for $Q$ to denote $|Q|={\bf sgn}(Q)\cdot Q$.

For a matrix $A\in\mathbb{R}^{n\times n}$, ${\bf span}(A)$ denotes
the column space of $A$, which is the linear span of its column vectors.
Given a set of vectors $\mathcal{B}=\{v_{1},v_{2},...,v_{r}\},v_{i}\in\mathbb{R}^{n},i\in\underline{r}$,
we define the matrix $B=[v_{1},v_{2},...,v_{r}]$ that stacks these
vectors, assuming the vectors are well ordered. For a matrix $A\in\mathbb{R}^{n\times n}$
with $rank(A)<n$, let $\mathcal{B}_{A}=\{\zeta_{1},\zeta_{2},...,\zeta_{m}\}$
represent the set of its eigenvectors corresponding to the eigenvalue
zero. This gives ${\bf null}(A)={\bf span}(\mathcal{B}_{A})$ along
with $B_{A}=[\zeta_{1},\zeta_{2},...,\zeta_{m}]$. The notation ${\bf diag}\{\cdot\}$
represents a diagonal matrix whose diagonal elements are taken from
the sequence in $\{\cdot\}$. Meanwhile, ${\bf blkdiag}\{\cdot\}$
denotes a block diagonal matrix whose off-diagonal blocks are zero
matrices, and the matrices in $\{\cdot\}$ form the diagonal blocks.
For subsets $\mathcal{X},\mathcal{Y}$ of $\mathbb{R}^{n}$, the Minkowski
sum is expressed as $\mathcal{X}+\mathcal{Y}=\{x+y|x\in\mathcal{X},y\in\mathcal{Y}\}$.
Additionally, we define the set $k\mathcal{X}\coloneqq\{kx|x\in\mathcal{X}\}$
for $k\in\mathbb{R}$.

\subsection{Graph Theory}

(\textbf{Graph}) Define the matrix-weighted graph/network $\mathcal{G}$
as a triplet $\mathcal{G}=(\mathcal{V},\mathcal{E},\mathcal{A})$,
where $\mathcal{V}$ has order $N$. The simple graph, which is the
scope of this work, is a graph where the connections are undirected
and are without self-loop or multiple edges. A node in $\mathcal{V}$
is also referred to as an agent $\tau$, where $\tau$ is a bijection
on $\underline{N}$ and $\tau(i)$ will be written as $\tau_{i}$.
An edge is then the pair $e_{ij}\coloneqq(\tau_{i},\tau_{j}),\tau_{i},\tau_{j}\in\mathcal{V}$.
The function $\mathcal{W}:\mathcal{E}\rightarrow\mathcal{A}$ maps
the edge to its associated weight matrix, i.e., $\mathcal{W}(e_{ij})=\mathcal{W}((\tau_{i},\tau_{j}))=A_{ij}\in\mathcal{A},\forall e_{ij}\in\mathcal{E}$.

\noindent (\textbf{Laplacian}) In this paper, the weight matrices
$A_{ij}\in\mathbb{R}^{d\times d}$ are considered to be real symmetric,
which furthermore have $|A_{ij}|\succeq0$ or $|A_{ij}|\succ0$ if
$(\tau_{i},\tau_{j})\in\mathcal{E}$, and $A_{ij}=0$ otherwise, for
all $\tau_{i},\tau_{j}\in\mathcal{V}$. When $|A_{ij}|\succeq0$ we
refer to the edge $e_{ij}$ as a semidefinite edge, while to $e_{ij}$
with $|A_{ij}|\succ0$, a definite edge. The adjacency matrix for
a matrix-weighted graph $A=[A_{ij}]\in\mathbb{R}^{dN\times dN}$ is
then a block matrix such that the block on the $\tau_{i}$-th row
and the $\tau_{j}$-th column is $A_{ij}$. Let $\mathcal{N}_{\tau_{i}}=\left\{ \tau_{j}\in\mathcal{V}\,|\,(\tau_{i},\tau_{j})\in\mathcal{E}\right\} $
be the neighbor set of an agent $\tau_{i}\in\mathcal{V}$. We use
$C=\text{{\bf blkdiag}}\left\{ C_{1},C_{2},\cdots,C_{N}\right\} \in\mathbb{R}^{dN}$
to represent the matrix-weighted degree matrix of a graph where $C_{\tau_{i}}=\sum_{\tau_{j}\in\mathcal{N}_{\tau_{i}}}|A_{ij}|\in\mathbb{R}^{d\times d}$.
The matrix-valued Laplacian matrix of a matrix-weighted graph is defined
as $L(\mathcal{G})=C-A$, which is real symmetric. A gauge matrix
in this paper is defined by the diagonal block matrix $D=S\otimes I_{d}$,
where $S={\bf diag}\{\sigma_{1},\sigma_{2},...,\sigma_{N}\},$ $\sigma_{1}=1,\sigma_{i}=\pm1$,
$i\in\{2,3,...,N\}.$ 

\noindent (\textbf{Path}) On a matrix-weighted graph $\mathcal{G}=(\mathcal{V},\mathcal{E},\mathcal{A})$,
a path $\mathcal{P}$ is defined as a sequence of edges as $\{(\tau_{1},\tau_{2}),(\tau_{2},\tau_{3}),\ldots,(\tau_{p-1},\tau_{p})\}$
where $(\tau_{i},\tau_{i+1})\in\mathcal{E},i\in\underline{p-1},$
$\tau_{1},\tau_{2},...,\tau_{p}\in\mathcal{V}$ are all distinct and
$\tau_{1},\tau_{p}$ are called endpoints. Two paths are node-independent
if the nodes they traverse have none in common. The sign of a path
is defined as ${\bf sgn}(\mathcal{P}):=\prod_{i=1}^{\text{card}(\mathcal{P})}{\bf sgn}(A_{i,i+1})$.
The null space of the path is defined as ${\bf null}(\mathcal{P}):={\bf null}(A_{12})+{\bf null}(A_{23})+...+{\bf null}(A_{p-1,p})$,
where $+$ is the Minkowski sum. We say $\mathcal{P}$ is a definite
path if $\forall(\tau_{i},\tau_{i+1})\in\mathcal{P},|A_{i,i+1}|\succ0$,
and $\mathcal{P}$ is a semidefinite path if $\forall(\tau_{i},\tau_{i+1})\in\mathcal{P},|A_{i,i+1}|\succeq0$
and $A_{i,i+1}\neq0$. 

\noindent (\textbf{Tree. Cycle}.) A positive-negative tree $\mathcal{T}$
on a matrix-weighted graph $\mathcal{G}$ is a tree such that $\forall(\tau_{i},\tau_{j})\in\mathcal{T}$,
it satisfies that $|A_{ij}|\succ0$. A positive-negative spanning
tree of $\mathcal{G}$ is a positive-negative tree traversing all
nodes in $\mathcal{G}$. A cycle $\mathcal{C}$ of $\mathcal{G}$
is a path that has the same node as endpoints, i.e., $\mathcal{C}=\{(\tau_{1},\tau_{2}),(\tau_{2},\tau_{3}),\ldots,(\tau_{p-1},\tau_{1})\}$. 
\begin{defn}[\textbf{Weakly Connected}]
A matrix-weighted network $\mathcal{G}=(\mathcal{V},\mathcal{E},\mathcal{A})$
is called weakly connected if $\mathcal{G}$ does not have a positive-negative
spanning tree.
\end{defn}
\begin{defn}[\textbf{Structural Balance}]
A matrix-weighted network $\mathcal{G}=(\mathcal{V},\mathcal{E},\mathcal{A})$
is $(\mathcal{V}_{1},\mathcal{V}_{2})-$structurally balanced if there
exists a bipartition of nodes $\mathcal{V}=\mathcal{V}_{1}\cup\mathcal{V}_{2},\mathcal{V}_{1}\cap\mathcal{V}_{2}=\emptyset$,
such that the matrix-valued weight between any two nodes within each
subset is positive (semi-)definite, but negative (semi-)definite for
edges connecting nodes of different subsets. A matrix-weighted network
is structurally imbalanced if it is not structurally balanced.
\end{defn}
\begin{defn}[\textbf{Nontrivial Balancing Set}]
\noindent \label{def:NBS}Let $\mathcal{G}=(\mathcal{V},\mathcal{E},\mathcal{A})$
be a matrix-weighted network and $(\mathcal{V}_{1},\mathcal{V}_{2})$
be a partition of $\mathcal{V}$. A Nontrivial Balancing Set (NBS)
$\mathcal{E}^{nb}(\mathcal{V}_{1},\mathcal{V}_{2})\subset\mathcal{E}$
defines a set of edges of $\mathcal{G}$ such that, (I) the negation
of the signs of their weight matrices renders $\mathcal{G}$ a $(\mathcal{V}_{1},\mathcal{V}_{2})$\textendash structurally
balanced graph, and (II) the null spaces of their weight matrices
have intersection other than $\{{\bf 0}\}$, which is denoted as ${\bf null}(\mathcal{E}^{nb}(\mathcal{V}_{1},\mathcal{V}_{2}))$. 
\end{defn}
\begin{rem}
\noindent When $\mathcal{E}^{nb}(\mathcal{V}_{1},\mathcal{V}_{2})=\emptyset$
because $\mathcal{G}$ is itself $(\mathcal{V}_{1},\mathcal{V}_{2})$\textendash structurally
balanced, we define in this case ${\bf null}(\mathcal{E}^{nb}(\mathcal{V}_{1},\mathcal{V}_{2}))={\bf null}({\bf 0}_{d\times d})=\mathbb{R}^{d\times d}$
suppose $\mathcal{A}\subset\mathbb{R}^{d\times d}$. An NBS is unique
in $\mathcal{G}$ if for any other partition $(\mathcal{V}_{1}',\mathcal{V}_{2}')$
of $\mathcal{V}$, such set of edges does not exist.
\end{rem}

\subsection{Problem Statement}

Consider a multi-agent network of $N\in\mathbb{Z}_{+}$ agents. The
states of each agent $\tau_{i}\in\mathcal{V}$ is denoted by $x_{\tau_{i}}(t)=\mathbb{R}^{d}$
where $d\in\mathbb{Z}_{+}$. The interaction protocol for the matrix-weighted
network reads

\begin{equation}
\dot{x}_{\tau_{i}}(t)=-\sum_{\tau_{j}\in\mathcal{N}_{\tau_{i}}}|A_{ij}|(x_{\tau_{i}}(t)-\text{{\bf sgn}}(A_{ij})x_{\tau_{j}}(t)),\tau_{i}\in\mathcal{V},\label{eq:protocol}
\end{equation}

\noindent where $A_{ij}\in\mathbb{R}^{d\times d}$ denotes the weight
matrix on edge $(\tau_{i},\tau_{j})$. The collective dynamics of
the multi-agent network \eqref{eq:protocol} can be characterized
by

\begin{equation}
\dot{x}(t)=-Lx(t),\label{eq:overall-dynamics}
\end{equation}
where $x(t)=[x_{1}^{T}(t),x_{2}^{T}(t),\ldots,x_{N}^{T}(t)]^{T}\in\mathbb{R}^{dN}$
and $L$ is the matrix-valued graph Laplacian. The control objective
for the system \eqref{eq:protocol} is related to the following definition. 
\begin{defn}[\textbf{Bipartite Consensus}]
\label{def:bi-consensus}The system \eqref{eq:protocol} achieves
bipartite consensus if there exist a partition $\mathcal{V}=\mathcal{V}_{1}\dot{\cup}\mathcal{V}_{2}$
and $x_{c}\in\mathbb{R}^{d},x_{c}\neq{\bf 0}$, such that, for all
initial values other than $x(0)={\bf 0},$ there is $\lim_{t\rightarrow\infty}x_{\tau_{i}}(t)=x_{c},\forall\tau_{i}\in\mathcal{V}_{1},$
and $\lim_{t\rightarrow\infty}x_{\tau_{i}}(t)=-x_{c}$, $\forall\tau_{i}\in\mathcal{V}_{2}$.
\end{defn}
The control objective is to determine the configuration of the network
\eqref{eq:protocol}, in terms of its structure and weight matrices,
that is sufficient to guarantee an asymptotic bipartite consensus
solution, given that the network is weakly connected.

\subsection{Network Properties}

The following results establish fundamental properties of network
\eqref{eq:protocol} that are essential for proving our main results.
\begin{lem}
\label{lem:convergence}For matrix-weighted network $\mathcal{G}$
in \eqref{eq:protocol} whose Laplacian has ${\bf null}(L)={\bf span}\{\eta_{1},\eta_{2},...,\eta_{m}\}$,
where $\eta_{1},\eta_{2},...,\eta_{m}$ are the orthonormal eigenvectors,
$\underset{t\rightarrow\infty}{\lim}x(t)$ exists and $\underset{t\rightarrow\infty}{\lim}x(t)=x^{*}=\sum_{i=1}^{m}(\eta_{i}^{T}x(t_{0}))\eta_{i}$.
\end{lem}
\begin{IEEEproof}
According to \cite[Lemma 1]{wang2022characterizing}, the matrix-valued
Laplacian satisfies 
\[
L=H^{T}{\bf blkdiag}\{|A_{k}|\}H
\]
where $H$ is the signed incidence matrix; since ${\bf blkdiag}\{|A_{k}|\}$
is positive semi-definite, $\forall x\in\mathbb{R}^{n}$ there is
$x^{T}Lx\geqslant0$. The fact that the Laplacian matrix is real and
positive semidefinite guarantees that it is diagonalizable, implying
all Jordan blocks for eigenvalue zero having dimension one. This ensures
the system's convergence to the equilibrium $x^{*}$ that satisfies
$Lx^{*}={\bf 0}$. As for the value of $x^{*}$, consider the system
response $x(t)=e^{-L(t-t_{0})}x(t_{0})=Qe^{J(t-t_{0})}Q^{T}x(t_{0})$
where $Q=[\eta_{n},...,\eta_{m+1},\eta_{m},...,\eta_{1}]$ is the
orthogonal matrix that satisfies $L=Q{\bf diag}\{\lambda_{n},...,\lambda_{m+1},0,...,0\}Q^{T}$,
$\lambda_{n}\geqslant\lambda_{n-1}\geqslant...\geqslant\lambda_{m+1}>0$.
Therefore we have $e^{J(t-t_{0})}={\bf diag}\{e^{-\lambda_{n}(t-t_{0})},...,e^{-\lambda_{m+1}(t-t_{0})},1,...,1\}$,
and $\underset{t\rightarrow\infty}{\lim}x(t)=\sum_{i=1}^{m}(\eta_{i}^{T}x(t_{0}))\eta_{i}$.
\end{IEEEproof}
\begin{prop}
\emph{(\cite[Theorem 1]{wang2022characterizing}) \label{pro:sol_Aij}}
Consider the matrix-weighted network $\mathcal{G}(\mathcal{V},\mathcal{E},\mathcal{A})$
in \eqref{eq:protocol} and its overall dynamics \eqref{eq:overall-dynamics},
$Lx={\bf 0}$ if and only if 
\begin{equation}
A_{ij}(x_{\tau_{i}}-{\bf sgn}(A_{ij})x_{\tau_{j}})={\bf 0},\forall(\tau_{i},\tau_{j})\in\mathcal{E}.\label{eq:sol_Aij}
\end{equation}
\end{prop}
\begin{rem}
Lemma \ref{lem:convergence} and Proposition \ref{pro:sol_Aij} show
that the agents of the matrix-weighted network converge to the null
space of the network's Laplacian, while to study the structure of
${\bf null}(\mathcal{L})$ is to investigate the solution to the set
of equations \eqref{eq:sol_Aij} for all the connections in the graph.
\end{rem}
The following proposition demonstrates that the bipartite consensus
solution corresponds to a specific structure of the Laplacian null
space.
\begin{prop}
\emph{\label{prop:BC_subspace}(}\textup{\cite[Theorem 1]{su2019bipartite}}\emph{)
}The matrix-weighted network \textup{\eqref{eq:overall-dynamics}}
achieves bipartite consensus if and only if there exists a gauge matrix
$D$ such that ${\bf null}(L)={\bf span}\left(D\left({\bf 1}_{n}\otimes\Psi\right)\right)$,
where $\Psi=\begin{bmatrix}\psi_{1} & \psi_{2} & \cdots & \psi_{s}\end{bmatrix}$,
$s=\text{dim}\left({\bf null}\left(L\right)\right)$, and $\psi_{i},i\in\underline{s}$,
are unit basis vectors.
\end{prop}
In \cite{wang2022characterizing}, another interpretation of bipartite
consensus on the matrix-weighted network is provided from the graph-theoretic
perspective instead of the purely algebraic, which relies on the concept
of the Nontrivial Balancing Set (see Definition \ref{def:NBS}). For
now, we denote ${\bf null}(\mathcal{E}^{nb}(\mathcal{V}_{1},\mathcal{V}_{2}))={\bf span}(B_{(\mathcal{V}_{1},\mathcal{V}_{2})})$
where the columns of the matrix $B_{(\mathcal{V}_{1},\mathcal{V}_{2})}\in\mathbb{R}^{d\times r}$
span the space ${\bf null}(\mathcal{E}^{nb}(\mathcal{V}_{1},\mathcal{V}_{2}))$
of dimension $r$. We see that the uniqueness of the NBS constitutes
a first necessary condition of the bipartite consensus achieved on
\eqref{eq:protocol} with general topology.
\begin{prop}
\emph{\label{pro:(W,Theorem-1)}(\cite[Theorem 1]{wang2022characterizing})}
Consider the matrix-weighted network \eqref{eq:protocol} and its
dynamics \eqref{eq:overall-dynamics}. There exists an NBS $\mathcal{E}^{nb}(\mathcal{V}_{1},\mathcal{V}_{2})$
in $\mathcal{G}(\mathcal{V},\mathcal{E},\mathcal{A})$ if and only
if ${\bf span}\left(D({\bf 1}_{N}\otimes B_{(\mathcal{V}_{1},\mathcal{V}_{2})})\right)\subset{\bf null}(L)$
where $D\in\mathbb{R}^{Nd\times Nd}$ is a gauge matrix, and ${\bf span}(B_{(\mathcal{V}_{1},\mathcal{V}_{2})})={\bf null}(\mathcal{E}^{nb}(\mathcal{V}_{1},\mathcal{V}_{2}))$;
that is, such a set of edges in the graph interchanges with a set
of vectors in the Laplacian null space.
\end{prop}
\begin{prop}
\emph{\label{pro:(W,Theorem-2)}(\cite[Theorem 2]{wang2022characterizing})}
If the agents of matrix-weighted network \eqref{eq:protocol} achieve
bipartite consensus between $\mathcal{V}_{1}$ and $\mathcal{V}_{2}$,
then there exists a unique NBS $\mathcal{E}^{nb}(\mathcal{V}_{1},\mathcal{V}_{2})$
in $\mathcal{G}(\mathcal{V},\mathcal{E},\mathcal{A})$ such that $x_{\tau_{i}}(\infty)\in{\bf null}(\mathcal{E}^{nb}(\mathcal{V}_{1},\mathcal{V}_{2}))$
for all $\tau_{i}\in\mathcal{V}$.
\end{prop}
\begin{prop}
\emph{\label{pro:(W,Theorem-3)}(\cite[Theorem 3]{wang2022characterizing})}
For the matrix-weighted network (\ref{eq:protocol}) with a positive-negative
spanning tree, the system admits bipartite consensus if and only if
there exists a unique NBS $\mathcal{E}^{nb}$ in the graph $\mathcal{G}$,
and it admits a trivial consensus when there is no such sets in the
graph. 
\end{prop}
Unless otherwise specified, the matrix-weighted networks studied in
this work are meant to be connected. As we point out in the following
Lemma, this assumption is de facto necessary for bipartite consensus.
\begin{lem}
A matrix-weighted network \eqref{eq:protocol} whose underlying graph
is not connected does not achieve bipartite consensus in the sense
of Definition \ref{def:bi-consensus}.
\end{lem}
\begin{IEEEproof}
Assume that bipartite consensus is achieved in the sense of Definition
\ref{def:bi-consensus} on network \eqref{eq:protocol} that is unconnected
between $\mathcal{V}_{1}$ and $\mathcal{V}_{2}$, $\mathcal{V}=\mathcal{V}_{1}\dot{\cup}\mathcal{V}_{2}$.
Then there exists a nontrivial balancing set $\mathcal{E}^{nb}$ in
the network such that, by changing the signs of the matrix weights
of the edges in $\mathcal{E}^{nb}$, the network $\mathcal{G}$ is
rendered structurally balanced; since structural balance is a global
property, the subgraphs based on $\mathcal{V}_{1}$ and $\mathcal{V}_{2}$
that are independent of each other are also rendered structurally
balanced. This suggests that the subgraphs $\mathcal{G}_{1}(\mathcal{V}_{1})$,
$\mathcal{G}_{2}(\mathcal{V}_{2})$ both have their own nontrivial
balancing sets, and according to Proposition \ref{pro:(W,Theorem-1)},
${\bf span}\left(D_{1}({\bf 1}_{|\mathcal{V}_{1}|}\otimes B_{1})\right)\subset{\bf null}(L_{1})$,
${\bf span}\left(D_{2}({\bf 1}_{|\mathcal{V}_{2}|}\otimes B_{2})\right)\subset{\bf null}(L_{2})$,
where $L_{1},L_{2}$ are the Laplacians of $\mathcal{G}_{1}(\mathcal{V}_{1})$,
$\mathcal{G}_{2}(\mathcal{V}_{2})$, while the Laplacian of $\mathcal{G}$
can be arranged as 
\[
L=\begin{bmatrix}L_{1} & O\\
O & L_{2}
\end{bmatrix}.
\]
Thus we know 
\[
{\bf span}\left(\begin{bmatrix}k_{1}D_{1}({\bf 1}_{|\mathcal{V}_{1}|}\otimes B_{1})\\
k_{2}D_{2}({\bf 1}_{|\mathcal{V}_{2}|}\otimes B_{2})
\end{bmatrix}\right)\subset{\bf null}(L)
\]
for $k_{1},k_{2}\in\mathbb{R},$ and bipartite consensus of the entire
system is not achieved in the sense of Definition \ref{def:bi-consensus}.
Therefore the assumption does not stand, and bipartite consensus must
be achieved with the precondition that the network $\mathcal{G}$
is connected.
\end{IEEEproof}

\section{Main Results}

This section provides sufficient conditions for the agents to achieve
bipartite consensus on a weakly connected network \eqref{eq:protocol}.
Generally, a matrix-weighted network can be divided into the strongly
connected parts \textendash{} subgraphs spanned by positive-negative
trees \textendash{} and the weak, semidefinite, connections between
them. This section introduces the concept of a ``continent'' for
the strongly connected part. It is then discussed how the network
can achieve bipartite consensus when several ``continents'' are
linked by weak connections known as the semidefinite paths. 

\subsection{Continents}

This subsection defines a strongly connected component of a matrix-weighted
graph, called a continent, and discusses its structural and dynamic
properties related to the NBS. 
\begin{defn}[\textbf{Continent}]
\noindent \label{def:continent}Given a matrix-weighted network $\mathcal{G}(\mathcal{V},\mathcal{E},\mathcal{A})$
whose matrices are of dimension $d\times d$, a maximal positive-negative
tree $\mathcal{T}_{m}(\mathcal{V}_{m},\mathcal{E}_{m},\mathcal{A}_{m})$
is a positive-negative tree such that $\forall\tau\in\mathcal{V}_{m},\forall\tau'\in\mathcal{V}\backslash\mathcal{V}_{m}$,
there does not exist a definite path with $\tau$ and $\tau'$ as
endpoints. A \emph{continent} of $\mathcal{T}_{m}$ is then the subgraph
$\mathcal{K}(\bar{\mathcal{V}},\bar{\mathcal{E}},\bar{\mathcal{A}})$
where $\bar{\mathcal{V}}=\mathcal{V}_{m}$, $\bar{\mathcal{E}}=\{(\tau_{i},\tau_{j})\mid\tau_{i},\tau_{j}\in\mathcal{V}_{m},\mathcal{W}((\tau_{i},\tau_{j}))\neq{\bf 0}_{d\times d}\}$,
and $\bar{\mathcal{A}}=\{\mathcal{W}((\tau_{i},\tau_{j}))\mid(\tau_{i},\tau_{j})\in\bar{\mathcal{E}}\}$.
\end{defn}
Since a continent $\mathcal{K}$ is essentially a subgraph of $\mathcal{G}$,
the concept of NBS is directly applied to $\mathcal{K}$. Here we
discuss the properties of a continent as only part of a matrix-weighted
graph in the presence of a unique NBS. 

Lemma \ref{lem:continent_NBS} states that structurally, a continent
$\mathcal{K}$ must have a unique NBS if the entire graph $\mathcal{G}$
has a unique NBS, and that the two NBSs have identical partitions
for the agents of the continent.
\begin{lem}
\label{lem:continent_NBS}Assume that $\mathcal{G}(\mathcal{V},\mathcal{E},\mathcal{A})$
is a connected matrix-weighted network with a unique NBS $\mathcal{E}^{nb}(\mathcal{V}_{1},\mathcal{V}_{2})$,
and $\mathcal{K}(\bar{\mathcal{V}},\bar{\mathcal{E}},\bar{\mathcal{A}})$
is a continent of $\mathcal{G}$, then $\mathcal{K}$ also has a unique
NBS $\mathcal{E}_{\mathcal{K}}^{nb}(\bar{\mathcal{V}}_{1},\bar{\mathcal{V}}_{2})$.
Moreover, the bipartitions $(\mathcal{V}_{1},\mathcal{V}_{2})$ and
$(\bar{\mathcal{V}}_{1},\bar{\mathcal{V}}_{2})$ are identical for
all $\tau_{p}\in\bar{\mathcal{V}}$.
\end{lem}
\begin{IEEEproof}
Since an NBS for $\mathcal{G}$ exists, that is, after the negation
of the signs of every $e\in\mathcal{E}^{nb}$, $\mathcal{G}$ becomes
structurally balanced, then any subgraph of $\mathcal{G}$ can also
turn structurally balanced by negating in a subset of $\mathcal{E}^{nb}$,
where the edges have weights with nontrivially intersecting null spaces.
Thus an NBS for a subgraph of $\mathcal{G}$ exists, and there is
an NBS $\mathcal{E}_{\mathcal{K}}^{nb}(\bar{\mathcal{V}}_{1},\bar{\mathcal{V}}_{2})$
for $\mathcal{K}$. By definition, the continent $\mathcal{K}$ has
a maximal positive-negative tree $\mathcal{T}_{m}(\mathcal{V}_{m},\mathcal{E}_{m},\mathcal{A}_{m})$.
We will show that $\mathcal{E}_{\mathcal{K}}^{nb}(\bar{\mathcal{V}}_{1},\bar{\mathcal{V}}_{2})$
is unique on $\mathcal{K}$, and the bipartition $(\bar{\mathcal{V}}_{1},\bar{\mathcal{V}}_{2})$
conforms that given by $\mathcal{T}_{m}$.

The bipartition $(\widetilde{\mathcal{V}}_{1},\widetilde{\mathcal{V}}_{2})$
of $\bar{\mathcal{V}}$ given by $\mathcal{T}_{m}$ is defined as
such: $\forall\tau_{p},\tau_{q}\in\bar{\mathcal{V}}$, if the path
on $\mathcal{T}_{m}$ between $\tau_{p}$ and $\tau_{q}$ is positive,
then $\tau_{p},\tau_{q}\in\widetilde{\mathcal{V}}_{1}$; otherwise,
$\tau_{p}\in\widetilde{\mathcal{V}}_{1},\tau_{q}\in\widetilde{\mathcal{V}}_{2}$.
Since any weight matrix of $\mathcal{A}_{m}$ is positive or negative
definite and $\mathcal{E}_{\mathcal{K}}^{nb}$ exists, $\mathcal{E}_{\mathcal{K}}^{nb}$
must include edges other than those of $\mathcal{E}_{m}$ to have
${\bf null}(\mathcal{E}_{\mathcal{K}}^{nb})\neq\{{\bf 0}\}$, thus
the bipartition of $\mathcal{T}_{m}$ is undisturbed, and $\bar{\mathcal{V}}_{1}=\widetilde{\mathcal{V}}_{1},\bar{\mathcal{V}}_{2}=\widetilde{\mathcal{V}}_{2}$.
On the other hand, if there exists $\mathcal{E}_{\mathcal{K}}^{nb}(\bar{\mathcal{V}}_{1}',\bar{\mathcal{V}}_{2}')$
that yields any other bipartition, at least one edge of $\mathcal{E}_{m}$
must be selected into $\mathcal{E}_{\mathcal{K}}^{nb}(\bar{\mathcal{V}}_{1}',\bar{\mathcal{V}}_{2}')$
for negation, then there will be the contradiction ${\bf null}(\mathcal{E}_{\mathcal{K}}^{nb}(\bar{\mathcal{V}}_{1}',\bar{\mathcal{V}}_{2}'))=\{{\bf 0}\}$.
For the same reason, the bipartition $\mathcal{E}^{nb}(\mathcal{V}_{1},\mathcal{V}_{2})$
made for $\tau_{p}\in\bar{\mathcal{V}}$ conforms that of $\mathcal{T}_{m}$. 
\end{IEEEproof}
\begin{cor}
For a matrix-weighted network $\mathcal{G}$ with a unique nontrivial
balancing set $\mathcal{E}^{nb}$ and continents $\mathcal{K}_{i}(\bar{\mathcal{V}}_{i},\bar{\mathcal{E}}_{i},\bar{\mathcal{A}}_{i}),i\in\underline{w}$,
if an edge $e\in\mathcal{E}_{\mathcal{K}_{i}}^{nb}$, then $e\in\mathcal{E}^{nb}$.
In addition, there is 
\begin{eqnarray*}
 &  & {\bf null}(\mathcal{E}^{nb})\\
 & = & \left(\bigcap_{i=1}^{w}{\bf span}(\mathcal{B}_{\mathcal{K}_{i}})\right)\cap\left(\bigcap_{e\in\mathcal{E}^{nb}\backslash\cup_{i=1}^{w}\mathcal{E}_{\mathcal{K}_{i}}^{nb}}{\bf span}(\mathcal{B}_{\mathcal{W}(e)})\right).
\end{eqnarray*}
\end{cor}
Dynamically, the steady-states of the agents of the continent $\mathcal{K}$
can also be characterized by its NBS, which involves the NBS's nontrivially
intersecting null space. Let us first introduce a basis representation
of some NBS-related null spaces that will come in handy for the mathematical
statements.
\begin{defn}[\textbf{Basis for Null Space}]
\noindent For a matrix-weighted network $\mathcal{G}$, if there
exists an NBS $\mathcal{E}^{nb}(\mathcal{V}_{1},\mathcal{V}_{2})$
with a nontrivially intersecting null space ${\bf null}(\mathcal{E}^{nb}(\mathcal{V}_{1},\mathcal{V}_{2}))$,
denote the basis of this null space as set $\mathcal{B}_{(\mathcal{V}_{1},\mathcal{V}_{2})}$,
which corresponds to a matrix $B_{(\mathcal{V}_{1},\mathcal{V}_{2})}$
whose columns are the basis vectors. Suppose $\mathcal{K}_{1},\mathcal{K}_{2},...,\mathcal{K}_{w}$
are continents of $\mathcal{G}$. For $\mathcal{K}_{i},i\in\underline{w}$,
if an NBS of the subgraph $\mathcal{K}_{i}$ exists, denote it as
$\mathcal{E}_{\mathcal{K}_{i}}^{nb}$ whose nontrivially intersecting
null space is ${\bf null}(\mathcal{E}_{\mathcal{K}_{i}}^{nb})$. Similarly,
denote the basis of ${\bf null}(\mathcal{E}_{\mathcal{K}_{i}}^{nb})$
as set $\mathcal{B}_{\mathcal{K}_{i}}$, whose corresponding matrix
is $B_{\mathcal{K}_{i}}$. 
\end{defn}
Combining Proposition \ref{pro:(W,Theorem-3)} and Lemma \ref{lem:continent_NBS},
we derive the following lemma about the steady-states of the agents
on the continent.
\begin{lem}
\label{lem:continent_convergence}Assume that $\mathcal{G}(\mathcal{V},\mathcal{E},\mathcal{A})$
is a connected matrix-weighted network with a unique NBS $\mathcal{E}^{nb}(\mathcal{V}_{1},\mathcal{V}_{2})$,
and $\mathcal{K}(\bar{\mathcal{V}},\bar{\mathcal{E}},\bar{\mathcal{A}})$
is a continent of $\mathcal{G}$, then the Laplacian matrix of the
subgraph $\mathcal{K}$ satisfies ${\bf null}(\mathcal{L}(\mathcal{K}))={\bf span}(\bar{D}({\bf 1}_{\text{card}\left(\bar{\mathcal{V}}\right)}\otimes B_{\mathcal{K}}))$,
where $\bar{D}\in\mathbb{R}^{Nd\times Nd}$ is a gauge matrix. Therefore,
within the entire system \eqref{eq:protocol}, the agents of the continent
$\mathcal{K}$ achieve bipartite consensus in the sense of Definition
\ref{def:bi-consensus} and the asymptotic value of agent $\tau_{p}\in\bar{\mathcal{V}}$
satisfies $x_{\tau_{p}}\in{\bf span}(\mathcal{B}_{\mathcal{K}})$.
\end{lem}
\begin{IEEEproof}
Proposition \ref{pro:(W,Theorem-3)} states that, a matrix-weighted
network with a positive-negative spanning tree admits bipartite consensus
if and only if there is a unique NBS in the graph, while to admit
the bipartite consensus solution is to have the null space of its
Laplacian ${\bf null}(\mathcal{L})$ as the bipartite consensus subspace
(Proposition \ref{prop:BC_subspace}). By Lemma \ref{lem:continent_NBS},
continent $\mathcal{K}$ has $\mathcal{T}_{m}$ along with a unique
NBS $\mathcal{E}_{\mathcal{K}}^{nb}$, thus $\mathcal{K}$ as a subgraph
meets ${\bf null}(\mathcal{L}(\mathcal{K}))={\bf span}(\bar{D}({\bf 1}_{|\bar{\mathcal{V}}|}\otimes B_{\mathcal{K}}))$
where $\bar{D}$ is a gauge matrix. This means that under the constraints
imposed by $\mathcal{K}$, the agents of $\mathcal{K}$ converge as
$x_{\tau_{p}}=\pm x_{\tau_{q}},x_{\tau_{p}}\in{\bf span}(B_{\mathcal{K}}),\tau_{p},\tau_{q}\in\bar{\mathcal{V}}$.
As we evaluate the whole system $\mathcal{G}$, what remains unchanged
is the bipartite solution $x_{\tau_{p}}=\pm x_{\tau_{q}},\tau_{p},\tau_{q}\in\bar{\mathcal{V}}$
obtained by solving eqn.\eqref{eq:sol_Aij} along $\mathcal{T}_{m}$.
Meanwhile, as we solve eqn.\eqref{eq:sol_Aij} on $(\tau_{p},\tau_{q}),\tau_{p}\in\bar{\mathcal{V}},\tau_{q}\notin\bar{\mathcal{V}}$,
additional constraints are put on $x_{\tau_{p}}$ thus its solution
space may reduce from ${\bf span}(\mathcal{B}_{\mathcal{K}})$; nevertheless,
the bipartition $x_{\tau_{p}}=\pm x_{\tau_{q}}$ for $\tau_{p},\tau_{q}\in\bar{\mathcal{V}}$
remains because of $\mathcal{T}_{m}$, and the statement $x_{\tau_{p}}\in{\bf span}(\mathcal{B}_{\mathcal{K}})$
stands. 
\end{IEEEproof}

\subsection{Continents Connected with Semidefinite Paths}

Now the semidefinite paths connecting the continents are considered
to determine the steady-state solution of the entire system. Without
loss of generality, the semidefinite paths discussed here do not contain
nodes from the continents, apart from their endpoints. Suppose the
matrix-weighted network \eqref{eq:protocol} has continents $\mathcal{K}_{1},...,\mathcal{K}_{\lambda}$.
It can be shown that for any pair of continents $(\mathcal{K}_{l},\mathcal{K}_{m})$
connected with semidefinite paths $\mathcal{P}_{1},...,\mathcal{P}_{\mu}$,
the asymptotic values of agents $\tau_{l}\in\mathcal{V}_{l},\tau_{m}\in\mathcal{V}_{m}$
satisfy the following equations:
\begin{equation}
\begin{array}{cll}
x_{\tau_{l}}-x_{\tau_{m}}\in\mathbf{null}(\mathcal{P}_{r}), & r\in\{1,...,\mu_{1}\}, & (\mathrm{I})\\
x_{\tau_{l}}+x_{\tau_{m}}\in\mathbf{null}(\mathcal{P}_{r}), & r\in\{\mu_{1}+1,...,\mu\}. & (\mathrm{II})
\end{array}\label{eq:continent_path}
\end{equation}
According to Lemma \ref{lem:continent_convergence}, $x_{\tau_{l}}$
represents the convergence value of all agents in $\mathcal{K}_{l}$,
provided there is a unique NBS. Therefore, the solution to eqn. \eqref{eq:continent_path}
determines whether agents of continents $\mathcal{K}_{1},...,\mathcal{K}_{\lambda}$
achieve bipartite consensus, leading to the following sufficient condition
in Lemma \ref{lem:bipartite-continent} to ensure bipartite convergence
between different continents.
\begin{lem}
\label{lem:bipartite-continent}Assume that a matrix-weighted network
\eqref{eq:protocol} with a unique NBS has continents $\mathcal{K}_{1},...,\mathcal{K}_{\lambda}$.
Then the agents of the continents $\mathcal{K}_{1},...,\mathcal{K}_{\lambda}$
achieve bipartite consensus if, for any $\mathcal{K}_{l},\mathcal{K}_{m}$
connected with semidefinite paths $\mathcal{P}_{1},...,\mathcal{P}_{\mu}$,
\begin{equation}
\begin{cases}
{\bf span}\left(\mathcal{B}_{\mathcal{K}_{l}}\cup\mathcal{B}_{\mathcal{K}_{m}}\right)\cap\left(\cap_{r=\alpha}^{\alpha'}{\bf null}\left(\mathcal{P}_{r}\right)\right)=\{{\bf 0}\}, & (\dot{\mathrm{I}})\\
{\bf span}\left(\mathcal{B}_{\mathcal{K}_{l}}\cup\mathcal{B}_{\mathcal{K}_{m}}\right)\cap\left(\cap_{r=\beta}^{\beta'}{\bf null}\left(\mathcal{P}_{r}\right)\right)\neq\{{\bf 0}\} & (\dot{\mathrm{II}})
\end{cases}\label{eq:bipartite_continent_suff}
\end{equation}
hold for $\alpha=1,\alpha'=\mu_{1},\beta=\mu_{1}+1,\beta'=\mu$ or
$\alpha=\mu_{1}+1,\alpha'=\mu,\beta=1,\beta'=\mu_{1}$.
\end{lem}
\begin{IEEEproof}
According to Lemma \ref{lem:continent_convergence}, since the network
\eqref{eq:protocol} has a unique NBS, agents of the same continent
achieve bipartite consensus, that is, $\forall\tau_{l},\tau_{l}'\in\mathcal{V}_{l}$,
there is $x_{\tau_{l}}(\infty)=\pm x_{\tau_{l}'}(\infty)$. Let $\mathcal{P}_{r}=\{(\tau_{1}^{r},\tau_{2}^{r}),(\tau_{2}^{r},\tau_{3}^{r}),...,(\tau_{\text{card}(\mathcal{P}_{r})}^{r},\tau_{\text{card}(\mathcal{P}_{r})+1}^{r})\}$,
$r\in\underline{\mu}$, where $\tau_{1}^{r}\in\mathcal{V}_{l}$, $\tau_{\text{card}(\mathcal{P}_{r})+1}^{r}\in\mathcal{V}_{m}$,
and $\{\tau_{2}^{r},\tau_{3}^{r},...,\tau_{\text{card}(\mathcal{P}_{r})}^{r}\}\cap\left(\bigcup_{u\in\underline{\lambda}}\mathcal{V}_{u}\right)=\emptyset$.
Then in light of Proposition \ref{pro:sol_Aij}, by cancelling out
$x_{\tau_{i}^{r}}-{\bf sgn}(A_{i,i+1}^{r})x_{\tau_{i+1}^{r}}\in{\bf null}(A_{i,i+1}^{r})$
along all the semidefinite paths between $\mathcal{K}_{l}$ and $\mathcal{K}_{m}$,
we show that $x_{\tau_{1}^{r}},x_{\tau_{\text{card}(\mathcal{P}_{r})+1}^{r}},r\in\underline{\mu}$,
are subject to constraints 
\[
\begin{array}{cll}
x_{\tau_{1}^{r}}-x_{\tau_{\text{card}(\mathcal{P}_{r})+1}^{r}}\in\mathbf{null}(\mathcal{P}_{r}), & r\in\{1,...,\mu_{1}'\}, & (\mathrm{I})\\
x_{\tau_{1}^{r}}+x_{\tau_{\text{card}(\mathcal{P}_{r})+1}^{r}}\in\mathbf{null}(\mathcal{P}_{r}), & r\in\{\mu_{1}'+1,...,\mu\}. & (\mathrm{II})
\end{array}
\]
This in turn means that $x_{\tau_{l}}$ and $x_{\tau_{m}}$ are subject
to constraints (I) and (II) for any $\tau_{l}\in\mathcal{V}_{l},\tau_{m}\in\mathcal{V}_{m}$.
It is then straightforward that all the agents of continents $\mathcal{K}_{1},...,\mathcal{K}_{\lambda}$
achieve bipartite consensus if for any $\mathcal{K}_{l},\mathcal{K}_{m},$
exactly one set of equations between $(\mathrm{I})$ and $(\mathrm{II})$
have a trivial solution, yielding $x_{\tau_{l}}=x_{\tau_{m}}$ or
$x_{\tau_{l}}=-x_{\tau_{m}}$, $\forall\tau_{l}\in\mathcal{V}_{l},\forall\tau_{m}\in\mathcal{V}_{m}$.

Now we prove the sufficiency statement. Due to Lemma \ref{lem:continent_convergence},
there is $x_{\tau_{l}}\pm x_{\tau_{m}}\in{\bf span}\left(\mathcal{B}_{\mathcal{K}_{l}}\cup\mathcal{B}_{\mathcal{K}_{m}}\right)$.
Suppose $\alpha=1,\alpha'=\mu_{1},\beta=\mu_{1}+1,\beta'=\mu$, then
$(\dot{\mathrm{I}})$ in eqn. \eqref{eq:bipartite_continent_suff}
guarantees that $x_{\tau_{l}}=x_{\tau_{m}}$, while $(\dot{\mathrm{II}})$
ensures that $x_{\tau_{l}}\neq{\bf 0}$, thus bipartite consensus
for agents on the continents is achieved. Suppose $\alpha=\mu_{1}+1,\alpha'=\mu,\beta=1,\beta'=\mu_{1}$,
then eqn. \eqref{eq:bipartite_continent_suff} guarantees $x_{\tau_{l}}=-x_{\tau_{m}}$
similarly, which completes the proof.
\end{IEEEproof}
\begin{rem}
\label{rem:necessary}Lemma \ref{lem:bipartite-continent} implies
that exactly one set of equations between (I) and (II) has a trivial
solution $\{{\bf 0}\}$ in order to achieve bipartite consensus.
\end{rem}
With the help of Lemma \ref{lem:bipartite-continent}, the main result
of this paper is captured in Theorem \ref{thm:semi-path}, with its
proof in Appendix.

\noindent \textbf{Assumption 1}: The matrix-weighted network $\mathcal{G}$
has a unique NBS $\mathcal{E}^{nb}(\mathcal{V}_{1},\mathcal{V}_{2})$.
\begin{thm}
\label{thm:semi-path}For a matrix-weighted network $\mathcal{G}$
with continents $\mathcal{K}_{1},...,\mathcal{K}_{\lambda}$, assume
Assumption 1 holds. If, for any $\mathcal{K}_{l},\mathcal{K}_{m}\in\{\mathcal{K}_{1},...,\mathcal{K}_{\lambda}\}$,
$\mathcal{G}$ meets the following conditions:

\emph{(1)} any semidefinite path between two arbitrary continents
$\mathcal{K}_{l}$ and $\mathcal{K}_{m}$ satisfies equation \eqref{eq:bipartite_continent_suff}. 

\emph{(2)} any two semidefinite paths connecting continents $\mathcal{K}_{l}$
and $\mathcal{K}_{m}$, excluding their endpoints, are node-independent;

\emph{(3)} for any semidefinite path $\mathring{\mathcal{P}}$ connecting
continents $\mathcal{K}_{l}$ and $\mathcal{K}_{m}$, its weight matrices
have linearly independent bases of null spaces, i.e., the columns
of $B_{A_{i,i+1}},(\tau_{i},\tau_{i+1})\in\mathring{\mathcal{P}},$
are linearly independent;

\emph{(4)} for a semidefinite path $\mathring{\mathcal{P}}$ connecting
continents $\mathcal{K}_{l}$ and $\mathcal{K}_{m}$ with $\mathring{\mathcal{P}}\cap\mathcal{E}^{nb}\neq\emptyset$,
the bases $B_{A_{i,i+1}},(\tau_{i},\tau_{i+1})\in\mathring{\mathcal{P}}\backslash\mathcal{E}^{nb}$,
and the basis of ${\bf span}(\mathcal{B}_{\mathcal{K}_{l}})\cap{\bf span}(\mathcal{B}_{\mathcal{K}_{m}})$
are linearly independent;

\noindent then under protocol (1) all the agents of $\mathcal{G}$
achieve bipartite consensus for almost all initial values.
\end{thm}
Theorem \ref{thm:semi-path} provides sufficient conditions to achieve
bipartite consensus in the presence of a weakly connected graph. While
these conditions are straightforward for engineering practitioners
to verify, it is important to emphasize that they are not overly conservative.

As shown in Proposition \ref{pro:(W,Theorem-2)}, Assumption 1 is
a necessary condition for achieving bipartite consensus. Once it is
satisfied, Conditions (1)\textendash (4) are sufficient to guarantee
bipartite consensus, though none of them is sufficient individually.
Within the four conditions, (2) plays a central role by imposing a
specific constraint on the network structure.

While Condition (1) provides sufficient condition that applies to
agents of the continents, it is easy to find examples where the violation
of (1) leads to steady-state distribution that is not bipartite consensus
(see Section \ref{sec:example}). When Assumption 1 and Conditions
(1) and (2) hold, Conditions (3) and (4) become both necessary and
sufficient to achieve bipartite consensus, thus do not further compromise
the generality of the theorem (see Claim \ref{claim:1}). This demonstrates
that the conditions presented in Theorem \ref{thm:semi-path} strike
a balance between necessity and practicality.

We now explain Condition (2) in more detail. Condition (2) requires
that any two semidefinite paths connecting two continents of the network
do not share any intermediate nodes. For instance, in Figure \ref{fig:eg1},
three semidefinite paths connect the two continents: $\mathcal{P}_{1}=\{(2,5)\},$$\mathcal{P}_{2}=\{(1,9),(9,4)\},$
$\mathcal{P}_{3}=\{(1,7),(7,8),(8,4)\}$. Although paths $\mathcal{P}_{2}$
and $\mathcal{P}_{3}$ have common endpoints (1 and 4), Condition
(2) requires that their intermediate nodes do not overlap, i.e. $\{7,8\}\cap\{9\}=\emptyset$.
However, if there exists an additional semidefinite edge $(9,5)$,
then the path $\{(1,9),(9,5)\}$ would share node 9 with path $\mathcal{P}_{2}$,
thereby violating Condition (2). 

Although Condition (2) narrows the focus to ensure sufficiency, it
is worth noting that alternative sufficient conditions can be derived
based on it to achieve bipartite consensus. For instance, when two
semidefinite paths share certain agents, provided these paths either
have no edge in the NBS (as described by eqn. \eqref{eq:gamma0_block})
or only one edge in the NBS (as described by eqn. \eqref{eq:gamma0'}),
our analysis remains valid. Under this new condition, bipartite consensus
can still be guaranteed when combined with Assumption 1 and Conditions
(1), (3) and (4). This further demonstrates that the conditions presented
in Theorem \ref{thm:semi-path} are not overly conservative.

\section{Simulation Example\label{sec:example}}

Numerical examples are presented in this section to validate Theorem
\ref{thm:semi-path} and compare it with existing results. 
\begin{example}
\label{exa:1}Network $\mathcal{G}_{1}$ in Figure \ref{fig:eg1}
consists of two continents, $\mathcal{K}_{1}$ of agents $\{1,2,3\}$,
$\mathcal{K}_{2}$ of agents $\{4,5,6\}$, and three connecting semidefinite
paths that are $\mathcal{P}_{1}=\{(2,5)\},$$\mathcal{P}_{2}=\{(1,9),(9,4)\},$
$\mathcal{P}_{3}=\{(1,7),(7,8),(4,8)\}$. The weight matrices are
determined with the following positive definite matrix $A_{def}$
and positive semidefinite matrices $A_{1},A_{2},...,A_{7}$, where 
\end{example}
\[
A_{def}=\begin{bmatrix}3 & 0 & 2 & 0\\
0 & 1 & 0 & 0\\
2 & 0 & 3 & 0\\
0 & 0 & 0 & 2
\end{bmatrix},A_{1}=\begin{bmatrix}2 & 2 & 0 & 0\\
2 & 2 & 0 & 0\\
0 & 0 & 1 & -1\\
0 & 0 & -1 & 1
\end{bmatrix},A_{2}=\begin{bmatrix}3 & 2 & 0 & 0\\
2 & 3 & 0 & 0\\
0 & 0 & 3 & -3\\
0 & 0 & -3 & 3
\end{bmatrix},
\]
\[
A_{3}=\begin{bmatrix}4 & 1 & 0 & 0\\
1 & 4 & 0 & 0\\
0 & 0 & 1 & -1\\
0 & 0 & -1 & 1
\end{bmatrix},A_{4}=\begin{bmatrix}2 & 1 & 0 & 0\\
1 & 2 & 0 & 0\\
0 & 0 & 1 & 1\\
0 & 0 & 1 & 1
\end{bmatrix},A_{5}=\begin{bmatrix}2 & -2 & 0 & 0\\
-2 & 2 & 0 & 0\\
0 & 0 & 5 & 1\\
0 & 0 & 1 & 5
\end{bmatrix},
\]
\[
A_{6}=\begin{bmatrix}3 & 0 & 1 & 0\\
0 & 0 & 0 & 0\\
1 & 0 & 2 & 0\\
0 & 0 & 0 & 2
\end{bmatrix}.
\]
We adopt $+A_{def}$($-A_{def}$) for the positive(negative) definite
matrices in network $\mathcal{G}_{1}$, and assign $\mathcal{W}((1,3))=-A_{1}$,
$\mathcal{W}((4,6))=A_{2}$, $\mathcal{W}((1,9))=-A_{3}$, $\mathcal{W}((2,5))=A_{4}$,
$\mathcal{W}((4,9))=\mathcal{W}((4,8))=-A_{4}$, $\mathcal{W}((1,7))=A_{5}$,
$\mathcal{W}((7,8))=A_{6}$. 

We now check for each condition in Theorem \ref{thm:semi-path}. 

\noindent \uline{Assumption 1} The network $\mathcal{G}_{1}$ is
structurally imbalanced due to the four negative cycles in the graph.
It is then observed that there exists a unique NBS $\mathcal{E}^{nb}(\mathcal{V}_{1},\mathcal{V}_{2})=\{(1,3),(4,6),(1,9)\}$
with $\mathcal{B}_{(\mathcal{V}_{1},\mathcal{V}_{2})}=\left\{ \begin{bmatrix}0 & 0 & 1 & 1\end{bmatrix}^{T}\right\} $,
whose negation of signs turns the graph structurally balanced between
$\mathcal{V}_{1}=\{1,3,6,7,8,9\},\mathcal{V}_{2}=\{2,4,5\}$.

\noindent \uline{Condition (1)} The two continents each has a semidefinite
edge $(1,3)$ or $(4,6)$ that constitutes its NBS, thus ${\bf span}(\mathcal{B}_{\mathcal{K}_{1}}\cup\mathcal{B}_{\mathcal{K}_{2}})={\bf span}(\left\{ \begin{bmatrix}1 & -1 & 0 & 0\end{bmatrix}^{T},\begin{bmatrix}0 & 0 & 1 & 1\end{bmatrix}^{T}\right\} )$.
Take the asymptotic values $x_{2}$ and $x_{5}$ as the representatives
of $\mathcal{K}_{1}$ and $\mathcal{K}_{2}$, we find the constraints
of semidefinite paths $\mathcal{P}_{1}$ and $\mathcal{P}_{3}$ to
be characterized by equation (I) in Lemma \ref{lem:bipartite-continent},
and $\mathcal{P}_{2}$, equation (II). As ${\bf span}(\mathcal{B}_{\mathcal{K}_{1}}\cup\mathcal{B}_{\mathcal{K}_{2}})\cap{\bf null}(\mathcal{P}_{1})\cap{\bf null}(\mathcal{P}_{3})=\{{\bf 0}\}$,
${\bf span}(\mathcal{B}_{\mathcal{K}_{1}}\cup\mathcal{B}_{\mathcal{K}_{2}})\cap{\bf null}(\mathcal{P}_{2})\neq\{{\bf 0}\}$,
equation \eqref{eq:bipartite_continent_suff} is satisfied.

\noindent \uline{Condition (2)} The semidefinite paths $\mathcal{P}_{1},\mathcal{P}_{2},\mathcal{P}_{3}$
are node independent except for the end points of $\mathcal{P}_{2}$
and $\mathcal{P}_{3}$. 

\noindent \uline{Condition (3)} It can be checked that the condition
holds.

\noindent \uline{Condition (4)} Semidefinite path $\mathcal{P}_{2}$
has $\mathcal{P}_{2}\cap\mathcal{E}^{nb}=(1,9)$, while the basis
of $\mathcal{W}((4,9))$ is linearly independent to that of ${\bf span}(\mathcal{B}_{\mathcal{K}_{1}})\cap{\bf span}(\mathcal{B}_{\mathcal{K}_{2}})$.

Thus network $\mathcal{G}_{1}$ satisfies all the assumption and conditions
of Theorem \ref{thm:semi-path}.

To verify whether the system dynamics converge to a bipartite consensus
solution, we introduce a metric $e_{b}(t)$ that measures the distance
of the system state from the bipartite consensus solution at time
$t$,
\[
e_{b}(t)=\sum_{k\in\mathcal{V}\backslash\{1\}}\begin{Vmatrix}x_{1}(t)-\sigma_{k}x_{k}(t)\end{Vmatrix}_{2}
\]
where $\sigma_{k}=\pm1$, and the sign pattern in sequence $\left(1,\sigma_{2},...,\sigma_{N}\right)$
indicates the bipartition to be expected in the bipartite consensus
solution.\footnote{By Proposition \ref{pro:(W,Theorem-2)}, if bipartite consensus is
achieved, the bipartition $(\mathcal{V}_{1},\mathcal{V}_{2})$ of
the agents is unique, and is indicated by the NBS $\mathcal{E}^{nb}(\mathcal{V}_{1},\mathcal{V}_{2})$
in the graph. Therefore one can obtain the sequence $\left(1,\sigma_{2},...,\sigma_{N}\right)$
once the graph is configured.} There is then $e_{b}(t)\rightarrow0$ if the agents achieve such
bipartite consensus, and $e_{b}(t)\rightarrow c>0$ otherwise.

\begin{figure}[t]
\begin{centering}
\begin{tikzpicture}[scale=1]
	\definecolor{dodgerblue}{RGB}{30,144,255} %135,206,250
	\definecolor{darkred}{RGB}{230,0,0}

    \tikzstyle{node1} = [circle,draw,minimum size=3mm]
	\tikzstyle{node2} = [circle,fill=white,minimum size=1mm]

	\node (n1) at (-1.8+5,1.8) [node1] {1};
	\node (n2) at (0+5,0) [node1] {2};
	\node (n3) at (-2.4+5,0) [node1] {3}; 
	\node (n4) at (3.6+5,1.8) [node1] {4}; 
	\node (n5) at (1.8+5,0) [node1] {5}; 
	\node (n6) at (4.2+5,0) [node1] {6}; 
	\node (n7) at (0+5,3) [node1] {7}; 
	\node (n8) at (1.8+5,3) [node1] {8}; 
	\node (n9) at (0.9+5,1.2) [node1] {9};

	\node (G1) at (-3.2+5,1.6) {$\mathcal{G}_{1}$};

	\tikzstyle{edge1} = [-, line width=1.4pt, color=darkred!90] % + def
	\tikzstyle{edge2} = [-, line width=1.4pt, color=dodgerblue] % - def
	\tikzstyle{edge3} = [-, line width=1.4pt, color=darkred!90, dashed] % + semi
	\tikzstyle{edge4} = [-, line width=1.4pt, color=dodgerblue, dashed] % - semi
	\draw[edge1]  (n4) -- (n5); 
	\draw[edge2]  (n5) -- (n6); 
	\draw[edge2]  (n2) -- (n3); 
	\draw[edge2]  (n1) -- (n2);
	\draw[edge3]  (n1) -- (n7);
	\draw[edge3]  (n7) -- (n8);
	\draw[edge4]  (n1) -- (n9);
	\draw[edge3]  (n6) -- (n4);
	\draw[edge3]  (n2) -- (n5);
	\draw[edge4]  (n4) -- (n8);
	\draw[edge4]  (n4) -- (n9);
	\draw[edge4]  (n1) -- (n3);

\end{tikzpicture}
\par\end{centering}
\caption{Network $\mathcal{G}_{1}$ and $\mathcal{G}_{3}$. Red solid (resp.,
dashed) lines represent edges weighted by positive definite (resp.,
semidefinite) matrices; blue solid (resp., dashed) lines represent
edges weighted by negative definite (resp., semidefinite) matrices.
Network $\mathcal{G}_{3}$ modifies $\mathcal{G}_{1}$ by letting
$\mathcal{W}((7,8))=A_{4}$, so that only Condition (3) does not hold. }

\label{fig:eg1}
\end{figure}
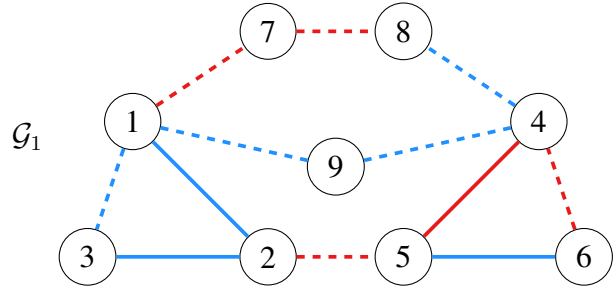

\begin{figure}[t]
\begin{centering}
\begin{tikzpicture}[scale=1]
	\definecolor{dodgerblue}{RGB}{30,144,255} %135,206,250
	\definecolor{darkred}{RGB}{230,0,0}

    \tikzstyle{node1} = [circle,draw,minimum size=3mm]
	\tikzstyle{node2} = [circle,fill=white,minimum size=1mm]

	\node (n1) at (-1.6,0) [node1] {1};
	\node (n2) at (-3.2,1.4) [node1] {2};
	\node (n3) at (-3.2,-1.4) [node1] {3}; 
	\node (n4) at (1.6,0) [node1] {4}; 
	\node (n5) at (3.2,-1.4) [node1] {5}; 
	\node (n6) at (3.2,1.4) [node1] {6}; 
	\node (n7) at (0,0) [node1] {7};

	\node (G2) at (-4.3,0) {$\mathcal{G}_{2}$};

	\tikzstyle{edge1} = [-, line width=1.4pt, color=darkred!90] % + def
	\tikzstyle{edge2} = [-, line width=1.4pt, color=dodgerblue] % - def
	\tikzstyle{edge3} = [-, line width=1.4pt, color=darkred!90, dashed] % + semi
	\tikzstyle{edge4} = [-, line width=1.4pt, color=dodgerblue, dashed] % - semi
	\draw[edge1]  (n1) -- (n3); 
	\draw[edge1]  (n1) -- (n2); 
	\draw[edge1]  (n4) -- (n5); 
	\draw[edge2]  (n4) -- (n6);
	\draw[edge3]  (n1) -- (n7);
	\draw[edge3]  (n4) -- (n7);
	\draw[edge3]  (n5) -- (n6);
	\draw[edge4]  (n3) -- (n2);

\end{tikzpicture}
\par\end{centering}
\caption{Network $\mathcal{G}_{2}$. Network $\mathcal{G}_{2}$ connects the
two continents with a single semidefinite path whose weight matrices
do not satisfy eqn. \eqref{eq:bipartite_continent_suff}, thereby
violating Condition (1). }

\label{fig:eg2}
\end{figure}
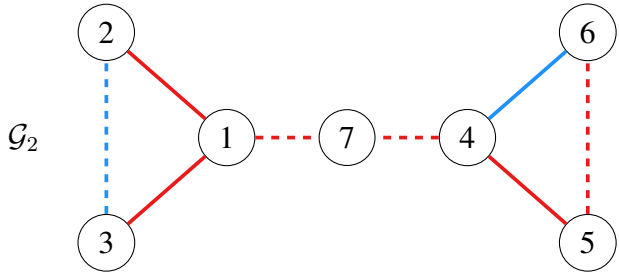

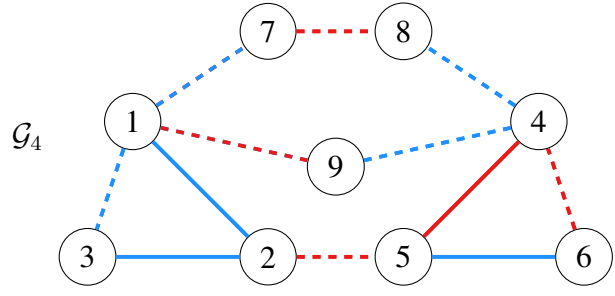
\begin{figure}[t]
\begin{centering}
\begin{tikzpicture}[scale=1]
	\definecolor{dodgerblue}{RGB}{30,144,255} %135,206,250
	\definecolor{darkred}{RGB}{230,0,0}

    \tikzstyle{node1} = [circle,draw,minimum size=3mm]
	\tikzstyle{node2} = [circle,fill=white,minimum size=1mm]

	\node (n1) at (-1.8+5,1.8) [node1] {1};
	\node (n2) at (0+5,0) [node1] {2};
	\node (n3) at (-2.4+5,0) [node1] {3}; 
	\node (n4) at (3.6+5,1.8) [node1] {4}; 
	\node (n5) at (1.8+5,0) [node1] {5}; 
	\node (n6) at (4.2+5,0) [node1] {6}; 
	\node (n7) at (0+5,3) [node1] {7}; 
	\node (n8) at (1.8+5,3) [node1] {8}; 
	\node (n9) at (0.9+5,1.2) [node1] {9};

	\node (G4) at (-3.2+5,1.6) {$\mathcal{G}_{4}$};

	\tikzstyle{edge1} = [-, line width=1.4pt, color=darkred!90] % + def
	\tikzstyle{edge2} = [-, line width=1.4pt, color=dodgerblue] % - def
	\tikzstyle{edge3} = [-, line width=1.4pt, color=darkred!90, dashed] % + semi
	\tikzstyle{edge4} = [-, line width=1.4pt, color=dodgerblue, dashed] % - semi
	\draw[edge1]  (n4) -- (n5); 
	\draw[edge2]  (n5) -- (n6); 
	\draw[edge2]  (n2) -- (n3); 
	\draw[edge2]  (n1) -- (n2);
	\draw[edge3]  (n1) -- (n7);
	\draw[edge3]  (n7) -- (n8);
	\draw[edge4]  (n1) -- (n9);
	\draw[edge3]  (n6) -- (n4);
	\draw[edge3]  (n2) -- (n5);
	\draw[edge4]  (n4) -- (n8);
	\draw[edge4]  (n4) -- (n9);
	\draw[edge4]  (n1) -- (n3);
	
	\draw[edge1]  (n4) -- (n5); 
	\draw[edge2]  (n5) -- (n6); 
	\draw[edge2]  (n2) -- (n3); 
	\draw[edge2]  (n1) -- (n2);
	\draw[edge4]  (n1) -- (n7);
	\draw[edge3]  (n7) -- (n8);
	\draw[edge3]  (n1) -- (n9);
	\draw[edge3]  (n6) -- (n4);
	\draw[edge3]  (n2) -- (n5);
	\draw[edge4]  (n4) -- (n8);
	\draw[edge4]  (n4) -- (n9);
	\draw[edge4]  (n1) -- (n3);

\end{tikzpicture}
\par\end{centering}
\caption{Network $\mathcal{G}_{4}$. Network $\mathcal{G}_{4}$ modifies $\mathcal{G}_{1}$
by letting $(1,7)\in\mathcal{E}^{nb}$, while the bases for the null
spaces of $\mathcal{W}((7,8)),\mathcal{W}((8,4)),\mathcal{W}((4,6))$
are linearly dependent, so that only Condition (4) does not hold.}

\label{fig:eg4}
\end{figure}

\begin{figure}[t]
\begin{centering}
\includegraphics[width=15cm]{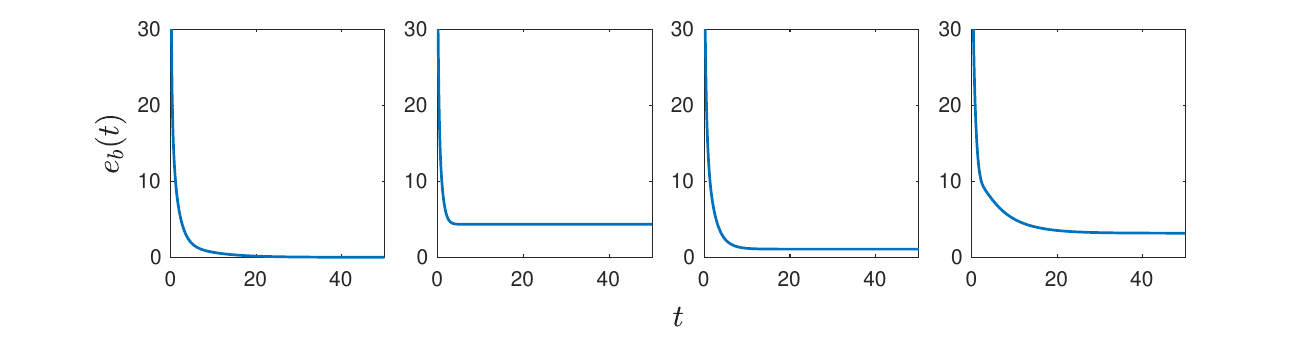}
\par\end{centering}
\caption{Metric $e_{b}(t)$ on $\mathcal{G}_{1}$, $\mathcal{G}_{2}$, $\mathcal{G}_{3}$,
and $\mathcal{G}_{4}$ (left to right). Bipartite consensus is achieved
on $\mathcal{G}_{1}$ but not on $\mathcal{G}_{2},\mathcal{G}_{3}$,
and $\mathcal{G}_{4}$.}

\label{fig:ebt}
\end{figure}

\begin{figure}[t]
\begin{centering}
\includegraphics[width=12cm]{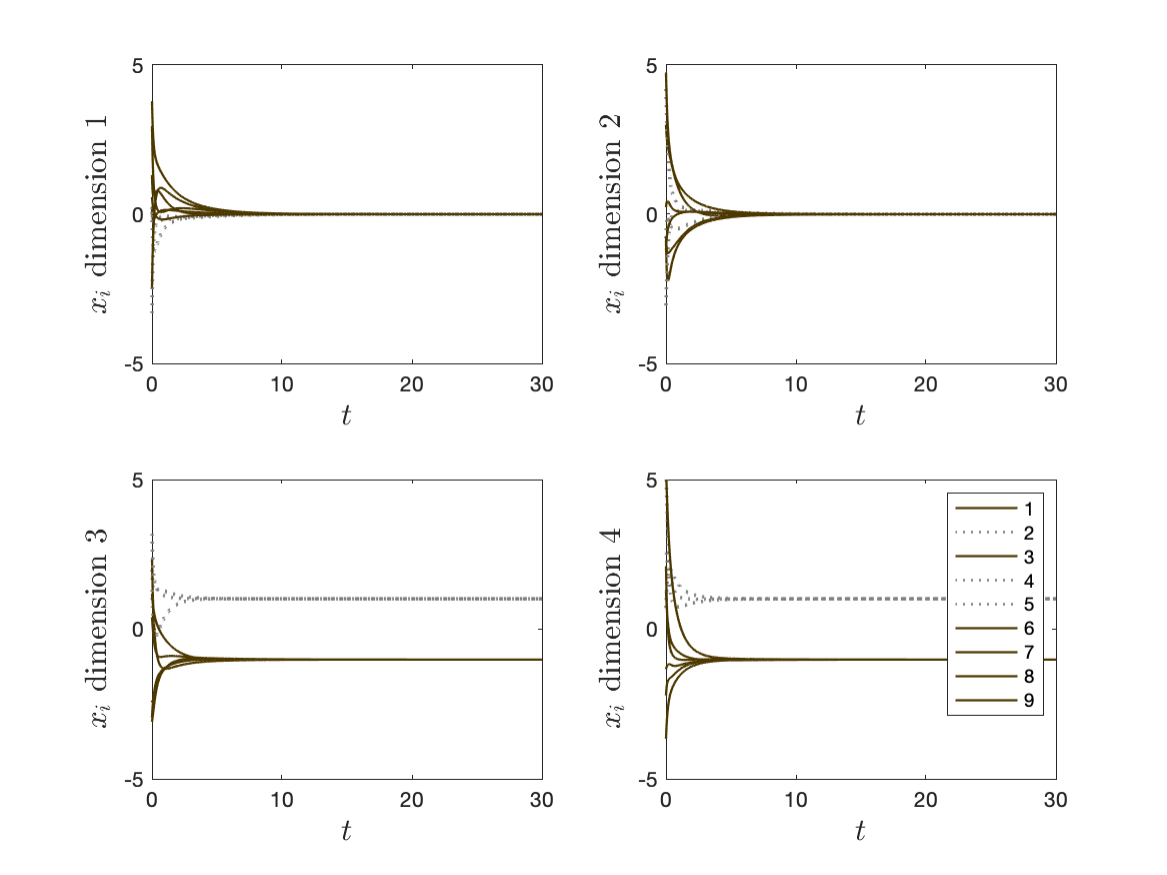}
\par\end{centering}
\caption{State trajectories of the nine agents of $\mathcal{G}_{1}$. Agents
of $\mathcal{V}_{1}=\{1,3,6,7,8,9\}$ and $\mathcal{V}_{2}=\{2,4,5\}$
are marked with solid lines and dashed lines respectively.}

\label{fig:eg1 state}
\end{figure}

The leftmost subfigure in Figure \ref{fig:ebt} demonstrates that
bipartite consensus is indeed achieved on network $\mathcal{G}_{1}$
as predicted by Theorem \ref{thm:semi-path}. The simulation result
in Figure \ref{fig:eg1 state} also shows that the bipartition occurs
between $\mathcal{V}_{1}=\{1,3,6,7,8,9\}$ and $\mathcal{V}_{2}=\{2,4,5\}$,
and the agents converge to ${\bf null}(\mathcal{E}^{nb})$. 

\begin{example}
In Example \ref{exa:1}, we can modify the structure or weight matrices
of $\mathcal{G}_{1}$ to construct a network that violates only one
of the three conditions. The other three subfigures of Figure \ref{fig:ebt}
depict the trajectory of $e_{b}(t)$ on modified networks $\mathcal{G}_{2},\mathcal{G}_{3}$
and $\mathcal{G}_{4}$ where Condition (1) or Condition (3) or Condition
(4) is not met, showing that bipartite consensus is not achieved.

Network $\mathcal{G}_{2}$, as depicted in Figure \ref{fig:eg2},
also has two continents on $\{1,2,3\}$ and $\{4,5,6\}$. The two
continents are connected by a semidefinite path $\mathcal{P}_{1}'=\{(1,7),(7,4)\}$.
The semidefinite edges have the following weights: $\mathcal{W}((2,3))=-A_{1}$,
$\mathcal{W}((5,6))=A_{2}$, $\mathcal{W}((1,7))=A_{7}$, $\mathcal{W}((4,7))=A_{8}$,
where
\[
A_{7}=\begin{bmatrix}2 & 2 & 0 & 0\\
2 & 2 & 0 & 0\\
0 & 0 & 3 & -1\\
0 & 0 & -1 & 3
\end{bmatrix},A_{8}=\begin{bmatrix}0 & 0 & 0 & 0\\
0 & 4 & 1 & 1\\
0 & 1 & 3 & 1\\
0 & 1 & 1 & 3
\end{bmatrix}.
\]
We then check each assumption and condition. Assumption 1 holds as
there exists a unique NBS $\mathcal{E}^{nb}(\mathcal{V}_{1},\mathcal{V}_{2})=\{(2,3),(5,6)\}$
between $\mathcal{V}_{1}=\{1,2,3,4,5,7\}$ and $\mathcal{V}_{2}=\{6\}$.
Condition (1) is violated because the semidefinite path is characterized
by equation (I) in Lemma \ref{lem:bipartite-continent}, while ${\bf null}\left(\mathcal{P}_{1}'\right)={\bf span}\left(\left\{ \begin{bmatrix}1 & 0 & 0 & 0\end{bmatrix}^{T},\begin{bmatrix}1 & -1 & 0 & 0\end{bmatrix}^{T}\right\} \right)$,
which intersects nontrivially with ${\bf span}(\mathcal{B}_{\mathcal{K}_{1}}\cup\mathcal{B}_{\mathcal{K}_{2}})={\bf span}\left(\left\{ \begin{bmatrix}1 & -1 & 0 & 0\end{bmatrix}^{T},\begin{bmatrix}0 & 0 & 1 & 1\end{bmatrix}^{T}\right\} \right)$.
Conditions (2) and (3) apparently hold, while Condition (4) does not
apply to this case as there is no edge on the semidefinite path belonging
to the NBS. Figure \ref{fig:ebt} shows that by failing Condition
(1) alone, network $\mathcal{G}_{2}$ does not achieve bipartite consensus.

Network $\mathcal{G}_{3}$ has the same structure as network $\mathcal{G}_{1}$,
but changes the weight matrix of edge $(7,8)$ to $\mathcal{W}((7,8))=A_{4}$.
This gives $\mathcal{W}((7,8))=-\mathcal{W}((4,8))$ and violates
Condition (3) on semidefinite path $\mathcal{P}_{2}$. One can check
that Assumption 1 and Conditions (1), (2), (4) are still satisfied.
Figure \ref{fig:ebt} indicates that bipartite consensus is not achieved
on network $\mathcal{G}_{3}$.

Network $\mathcal{G}_{4}$, as illustrated in Figure \ref{fig:eg4},
also has the same structure as network $\mathcal{G}_{1}$, but assigns
$\mathcal{W}((1,9))\succ0$ and $\mathcal{W}((1,7))\prec0$. The semidefinite
edges have the following weights: $\mathcal{W}((1,3))=-A_{1}$, $\mathcal{W}((4,6))=A_{1}$,
$\mathcal{W}((2,5))=A_{4}$, $\mathcal{W}((1,9))=A_{5}$, $\mathcal{W}((4,9))=-A_{6}$,
$\mathcal{W}((1,7))=-A_{3}$, $\mathcal{W}((7,8))=A_{6}$, $\mathcal{W}((4,8))=-A_{8}$.
The unique NBS of $\mathcal{G}_{3}$ is then $\mathcal{E}^{nb}(\mathcal{V}_{1},\mathcal{V}_{2})=\{(1,3),(4,6),(1,7)\}$,
with $\mathcal{V}_{1}=\{1,3,6,7,8,9\}$, $\mathcal{V}_{2}=\{2,4,5\}$.
One can check that Assumption 1 and Conditions (1), (2), (3) hold.
For Condition (4), $\mathcal{P}_{3}\cap\mathcal{E}^{nb}=\{(1,7)\}$,
but the bases of $\mathcal{W}((7,8))$ and $\mathcal{W}((4,8))$,
i.e., 
\[
\begin{bmatrix}0 & 1 & 0 & 0\end{bmatrix}^{T},\begin{bmatrix}1 & 0 & 0 & 0\end{bmatrix}^{T},
\]

\noindent and the basis of ${\bf span}(\mathcal{B}_{\mathcal{K}_{1}})\cap{\bf span}(\mathcal{B}_{\mathcal{K}_{2}})$,
i.e., 
\[
\begin{bmatrix}1 & -1 & 0 & 0\end{bmatrix}^{T},\begin{bmatrix}0 & 0 & 1 & 1\end{bmatrix}^{T},
\]
are linearly dependent. In this case, Figure \ref{fig:ebt} shows
that bipartite consensus is not achieved on network $\mathcal{G}_{4}$.
\end{example}
These additional examples substantiate that the sufficient conditions
presented in Theorem \ref{thm:semi-path} are not overly conservative.

\begin{rem}
Here we make a comparison between Theorem \ref{thm:semi-path} and
Theorem 4 in \cite{wang2022characterizing}. Much alike Theorem 3
in \cite{TRINH2018415}, Theorem 4 in \cite{wang2022characterizing}
proposes a sufficient condition to decide if an agent $\tau_{m}$
can achieve bipartite consensus with respect to an agent $\tau_{l}$
of a continent, or ``merge'' with it, given the various paths connecting
them. Therefore by applying this criterion to all agents outside the
continent, one can determine if an overall bipartite consensus can
be achieved. 

The main idea of both sufficiency theorems in \cite{wang2022characterizing}
and \cite{TRINH2018415} is that all the connecting paths $\mathcal{P}_{i}$
between $\tau_{m}$ and $\tau_{l}$ should have trivial intersection
of null spaces, i.e., $\bigcap_{i}{\bf null}(\mathcal{P}_{i})=\{{\bf 0}\}$,
so that the agents only yield $x_{\tau_{m}}=x_{\tau_{l}}$ or $x_{\tau_{m}}=-x_{\tau_{l}}$.
In the case of $\mathcal{G}_{1}$ where agents 4 and 5 naturally form
consensus as part of a continent, Theorem 4 of \cite{wang2022characterizing}
would take into account paths $\tilde{\mathcal{P}}_{1}=\{(2,5)\},$
$\tilde{\mathcal{P}}_{2}=\{(2,1),(1,9),(9,4)\},$ $\tilde{\mathcal{P}}_{3}=\{(2,1),(1,7),(7,8),(8,4)\}$,
$\tilde{\mathcal{P}}_{4}=\{(2,3),(3,1),(1,7),(7,8),(8,4)\}$ and $\tilde{\mathcal{P}}_{5}=\{(2,3),(3,1),(1,9),(9,4)\}$
to verify if agent 2 will merge into the cluster $\{4,5\}$. However,
it is obvious that $\bigcap_{i=1}^{5}{\bf null}(\tilde{\mathcal{P}}_{i})\neq\{{\bf 0}\}$
as edges $(2,5),(4,8)$, and $(4,9)$ all have the same null space
${\bf span}\left(\left\{ \begin{bmatrix}0 & 0 & 1 & -1\end{bmatrix}^{T}\right\} \right)$
to their weight matrices, despite that $\{2,4,5\}$ indeed forms a
cluster as indicated in Figure \ref{fig:eg1 state}. 

This comparison indicates that Theorem \ref{thm:semi-path} extends
previous work by relaxing the requirement for the null spaces of connecting
paths to intersect trivially. In fact, as ${\bf null}(\mathcal{P})$
is the Minkowski sum of the null spaces of its weight matrices, for
there to be $\bigcap_{i}{\bf null}(\mathcal{P}_{i})=\{{\bf 0}\}$,
the options for the semidefinite weight matrices on different paths
easily run out when $\tau_{m}$ has many paths to reach $\tau_{l}$.
Therefore, though not explicitly stated as conditions, strict topological
constraints on the semidefinite paths are implied in existing results.
Meanwhile, our result offers more flexibility regarding network structure,
such as the number and length of semidefinite paths, to achieve bipartite
consensus.
\end{rem}

\section{Concluding Remarks}

This work investigates sufficient conditions for achieving bipartite
consensus on weakly connected matrix-weighted networks. Recognizing
the difference between connections weighted by definite and semidefinite
matrices, we first introduce the concept of a ``continent'' \textendash{}
a subgraph of the network spanned by a positive-negative tree that
preserves its bipartite convergence regardless of the configuration
of the rest of the graph. Consequently, sufficient conditions for
achieving bipartite consensus on the entire network are presented,
allowing continents to be connected via the semidefinite paths. Finally,
the simulation results illustrate how to verify the proposed conditions
and demonstrate that the sufficient conditions are not overly conservative.

\appendix{}

Proof for Theorem \ref{thm:semi-path}.
\begin{IEEEproof}
Due to Assumption 1 which states that $\mathcal{G}$ has a unique
NBS, we have established in Lemma \ref{lem:continent_convergence}
that agents on each continent $\mathcal{K}_{l}$ converge bipartitely
to ${\bf null}(\mathcal{E}_{\mathcal{K}_{l}}^{nb})={\bf span}(\mathcal{B}_{\mathcal{K}_{l}})$.
Moreover, due to Condition (1) and Lemma \ref{lem:bipartite-continent},
bipartite consensus is achieved for agents from different continents
as well. The ambiguity rests with the agents on the semidefinite paths
that bridge the continents, for which we introduce Conditions (2)
\textendash{} (4). 

A formal statement of Condition (2) is given as follows. Assume there
are two semidefinite paths $\mathring{\mathcal{P}}_{1}=\{(\tau_{1}^{1},\tau_{2}^{1}),...,(\tau_{\rho}^{1},\tau_{\rho+1}^{1})\},\mathring{\mathcal{P}}_{2}=\{(\tau_{1}^{2},\tau_{2}^{2}),...,(\tau_{\rho'}^{2},\tau_{\rho'+1}^{2})\}$
between continents $\mathcal{K}_{l},\mathcal{K}_{m}$ where $\tau_{1}^{1},\tau_{1}^{2}\in\mathcal{K}_{l}$
and $\tau_{\rho+1}^{1},\tau_{\rho'+1}^{2}\in\mathcal{K}_{m}$, then
we allow $\tau_{1}^{1}=\tau_{1}^{2}$ and/or $\tau_{\rho+1}^{1}=\tau_{\rho'+1}^{2}$
without making it a necessity, but the rest of the path should be
node-independent which means $\{\tau_{2}^{1},\tau_{3}^{1},...,\tau_{\rho+1}^{1}\}\cap\{\tau_{2}^{2},\tau_{3}^{2},...,\tau_{\rho'+1}^{2}\}=\emptyset$.

We would first consider the scenario where there is only one semidefinite
path between two arbitrary continents, then extend to the general
situation of having multiple paths.

\noindent \uuline{Scenario 1}:

Assume that there is only one semidefinite path $\mathring{\mathcal{P}}=\{(\tau_{1},\tau_{2}),...,(\tau_{\rho},\tau_{\rho+1})\}$
between arbitrary continents $\mathcal{K}_{l},\mathcal{K}_{m}$. According
to Proposition \ref{pro:(W,Theorem-1)}, Assumption 1 guarantees that
$D({\bf 1}_{N}\otimes B_{(\mathcal{V}_{1},\mathcal{V}_{2})})$ is
included in the solution space for $\dot{x}={\bf 0}$ where $D$ specifies
the bipartition given by $\mathcal{E}^{nb}$; we will prove that the
columns of $D({\bf 1}_{N}\otimes B_{(\mathcal{V}_{1},\mathcal{V}_{2})})$
actually span the solution space by solving \eqref{eq:sol_Aij} for
$\mathring{\mathcal{P}}$. For the rest of the proof, we employ $\bar{D}={\bf diag}\{\sigma_{\hat{1}},\sigma_{\hat{2}},...,\sigma_{\hat{\rho}}\}\otimes I_{d}$
to indicate the bipartition of the agents that are on the semidefinite\textcolor{red}{{}
}path given by the NBS. The sign pattern of $\bar{D}$ conforms that
of $D$, i.e., $\sigma_{i}=\sigma_{\hat{i}}$ when $i$ and $\hat{i}$
denote the same agent. 

For $\mathring{\mathcal{P}}=\{(\tau_{1},\tau_{2}),...,(\tau_{\rho},\tau_{\rho+1})\}$,
since $\tau_{1}$ and $\tau_{\rho+1}$ converge bipartitely under
Assumption 1 and Condition (1), there is $x_{\tau_{1}}=s\cdot x_{\tau_{\rho+1}},s=\pm1,$
$x_{\tau_{1}},x_{\tau_{\rho+1}}\in{\bf span}(\mathcal{B}_{\mathcal{K}_{l}})\cap{\bf span}(\mathcal{B}_{\mathcal{K}_{m}})$.
We denote this solution for the end points as the null space of some
matrix $\bar{A}$, that is, 
\[
{\bf null}(\bar{A})={\bf span}(\mathcal{B}_{\mathcal{K}_{l}})\cap{\bf span}(\mathcal{B}_{\mathcal{K}_{m}}).
\]
To simplify the notation, we assign $\mathcal{W}((\tau_{i},\tau_{i+1}))=A_{i}$,
and its sign as ${\bf sgn}(A_{i})=s_{i}$. Then $s_{i}^{2}=1$, ${\bf sgn}(\mathring{\mathcal{P}})=s_{1}s_{2}\cdots s_{\rho}$.
Also, denote $\alpha_{i}=s_{1}s_{2}\cdots s_{i}$ for $i\in\underline{\rho}$.
Evaluate equation (\ref{eq:sol_Aij}) for the agents on $\mathring{\mathcal{P}}$,
there is 
\begin{align}
\Gamma_{0}x_{\mathring{\mathcal{P}}} & =\label{eq:gamma0}\\
 & \begin{bmatrix}\bar{A} & O & O & \cdots & O & O\\
I & O & O & \cdots & O & -sI\\
A_{1} & -s_{1}A_{1} & O & \cdots & O & O\\
O & A_{2} & -s_{2}A_{2} & \cdots & O & O\\
\vdots & \vdots & \vdots &  & \vdots & O\\
O & O & O & \cdots & A_{\rho} & -s_{\rho}A_{\rho}
\end{bmatrix}\begin{bmatrix}x_{\tau_{1}}\\
x_{\tau_{2}}\\
\vdots\\
x_{\tau_{\rho+1}}
\end{bmatrix}\nonumber \\
 & ={\bf 0}.\nonumber 
\end{align}
\uline{Case 1}: When $\mathring{\mathcal{P}}\cap\mathcal{E}^{nb}=\emptyset$,
the bipartition of $\mathring{\mathcal{P}}$ respects that of the
NBS, which gives $s={\bf sgn}(\mathring{\mathcal{P}})$. Then after
elementary column operations, $\Gamma_{0}$ can be turned into 
\begin{equation}
\begin{bmatrix}\bar{A} & O\\
O & R\\
O & Q
\end{bmatrix}=\begin{bmatrix}\bar{A} & O & O & \cdots & O\\
O & -\alpha_{1}I & -\alpha_{2}I & \cdots & -\alpha_{\rho}I\\
O & -s_{1}A_{1} & O & \cdots & O\\
O & O & -s_{2}A_{2} & \cdots & O\\
\vdots & \vdots & \vdots & \ddots & \vdots\\
O & O & O & \cdots & -s_{\rho}A_{\rho}
\end{bmatrix}\label{eq:gamma0_block}
\end{equation}
where $R=\begin{bmatrix}-\alpha_{1}I & -\alpha_{2}I & \cdots & -\alpha_{\rho}I\end{bmatrix}$,
$Q={\bf diag}\{-s_{1}A_{1},...,-s_{\rho}A_{\rho}\}$, $O$ stands
for the $d\times d$ zero matrix. It is easy to see that $\bar{D}({\bf 1}_{\rho}\otimes B_{\bar{A}})$
is a solution of \eqref{eq:gamma0}, and that if $\begin{bmatrix}R^{T} & Q^{T}\end{bmatrix}^{T}$
is of full rank, the nullity of $\Gamma_{0}$ in \eqref{eq:gamma0}
is equal to that of $\bar{A}$, then under the constraints of $\mathring{\mathcal{P}}$,
the nodes can only converge to ${\bf span}(\bar{D}({\bf 1}_{\rho}\otimes B_{\mathcal{K}_{l}}))\cap{\bf span}(\bar{D}({\bf 1}_{\rho}\otimes B_{\mathcal{K}_{m}}))$.
We then check if under Condition (3), there is 
\[
rank\left(\begin{bmatrix}R\\
Q
\end{bmatrix}\right)=\rho d.
\]
It is known that 
\begin{equation}
rank\left(\begin{bmatrix}R\\
Q
\end{bmatrix}\right)=rank(R)+rank(Q-QR^{\dagger}R)\label{eq:rank_block}
\end{equation}
where $R^{\dagger}$ is the Moore-Penrose inverse of $R$, and since
$R=\begin{bmatrix}-\alpha_{1} & -\alpha_{2} & \cdots & -\alpha_{\rho}\end{bmatrix}\otimes I$,
we have $R^{\dagger}=\frac{1}{\rho}\begin{bmatrix}-\alpha_{1} & -\alpha_{2} & \cdots & -\alpha_{\rho}\end{bmatrix}^{T}\otimes I$
and $rank(R)=d$. Therefore if $\bar{\Gamma}=Q(I-R^{\dagger}R)$ has
$rank(\bar{\Gamma})=(\rho-1)d$, we have the desired conclusion. Notice
that 
\[
\bar{\Gamma}=\begin{bmatrix}-s_{1}(1-\frac{1}{\rho})A_{1} & \frac{1}{\rho}\alpha_{2}A_{1} & \cdots & \frac{1}{\rho}\alpha_{\rho}A_{1}\\
\frac{1}{\rho}\alpha_{1}\alpha_{1}A_{2} & -s_{2}(1-\frac{1}{\rho})A_{2} & \cdots & \frac{1}{\rho}\alpha_{1}\alpha_{\rho}A_{2}\\
\vdots & \vdots &  & \vdots\\
\frac{1}{\rho}\alpha_{\rho-1}\alpha_{1}A_{\rho} & \frac{1}{\rho}\alpha_{\rho-1}\alpha_{2}A_{\rho} & \cdots & -s_{\rho}(1-\frac{1}{\rho})A_{\rho}
\end{bmatrix},
\]
let $\alpha=\begin{bmatrix}\alpha_{1} & \alpha_{2} & \cdots & \alpha_{\rho}\end{bmatrix}^{T}$
and $\varphi=\begin{bmatrix}\varphi_{1}^{T} & \varphi_{2}^{T} & \cdots & \varphi_{\rho}^{T}\end{bmatrix}^{T}$,
one can check that $\bar{\Gamma}\varphi={\bf 0}$ for $\varphi_{i}=\alpha_{i}v,i\in\underline{\rho},v\in\mathbb{R}^{d}$,
take the $\rho$-th row block for example, 
\begin{align*}
 & A_{\rho}\left(\frac{1}{\rho}\alpha_{\rho-1}v+\frac{1}{\rho}\alpha_{\rho-1}v+...+\frac{1}{\rho}s_{\rho}\alpha_{\rho}v-s_{\rho}\alpha_{\rho}v\right)\\
= & \alpha_{\rho-1}A_{\rho}v(\frac{1}{\rho}+\frac{1}{\rho}+...+\frac{1}{\rho}-1)\\
= & {\bf 0},
\end{align*}
which means ${\bf span}(\alpha\otimes I_{d})\subseteq{\bf null}(\bar{\Gamma})$,
therefore $\text{nullity}(\bar{\Gamma})\geqslant d$. Suppose there
exists $\psi=\begin{bmatrix}\psi_{1}^{T} & \psi_{2}^{T} & \cdots & \psi_{\rho}^{T}\end{bmatrix}^{T}$
such that $\psi\notin{\bf span}(\alpha\otimes I_{d})$ and $\psi\in{\bf null}(\bar{\Gamma}).$
Equation $\bar{\Gamma}\psi={\bf 0}$ is written as 
\[
\begin{array}{c}
s_{1}A_{1}(-\psi_{1}+\frac{\alpha_{1}}{\rho}\beta)={\bf 0}\\
s_{2}A_{2}(-\psi_{2}+\frac{\alpha_{2}}{\rho}\beta)={\bf 0}\\
\vdots\\
s_{\rho}A_{\rho}(-\psi_{1}+\frac{\alpha_{\rho}}{\rho}\beta)={\bf 0}
\end{array}
\]
where $\beta=\alpha_{1}\psi_{1}+\alpha_{2}\psi_{2}+...+\alpha_{\rho}\psi_{\rho}$,
which is equivalent to have $-\psi_{i}+\frac{\alpha_{i}}{\rho}\beta=\omega_{i}\in{\bf null}(A_{i})$
for $i\in\underline{\rho}$, meanwhile $\psi\notin{\bf span}(\alpha\otimes I_{d})$
requires that $\exists\omega_{i}\neq{\bf 0}$ for $i\in\underline{\rho}$.
However, the summation 
\[
\alpha_{1}\omega_{1}+\alpha_{2}\omega_{2}+...+\alpha_{\rho}\omega_{\rho}=-(\alpha_{1}\psi_{1}+\alpha_{2}\psi_{2}+...+\alpha_{\rho}\psi_{\rho})+\beta={\bf 0}
\]
implies that the nonzero $\omega_{i}$'s are linearly dependent which
contradicts Condition (3). We are then left with ${\bf null}(\bar{\Gamma})={\bf span}(\alpha\otimes I_{d})$
and $rank(\bar{\Gamma})=(\rho-1)d$. Therefore ${\bf null}(\Gamma_{0})={\bf span}(\bar{D}({\bf 1}_{\rho}\otimes B_{\bar{A}}))$
and the asymptotic values of the agents on $\mathring{\mathcal{P}}$
are in ${\bf span}(\mathcal{B}_{\mathcal{K}_{l}})\cap{\bf span}(\mathcal{B}_{\mathcal{K}_{m}})$.

\noindent \uline{Case 2}: If $\mathring{\mathcal{P}}\cap\mathcal{E}^{nb}\neq\emptyset$,
it is easy to see that $\mathring{\mathcal{P}}$ has exactly one edge
in the NBS because it is node-independent. In this case, the bipartition
of $\mathring{\mathcal{P}}$ does not respect that of the NBS and
we have $s=-{\bf sgn}(\mathring{\mathcal{P}})$, $\Gamma_{0}$ of
\eqref{eq:gamma0} instead becomes equal in rank to
\begin{equation}
\bar{\Gamma}_{0}=\begin{bmatrix}2I & \alpha_{1}I & \alpha_{2}I & \cdots & \alpha_{\rho-1}I & -\alpha_{\rho}I\\
\bar{A} & O & O & \cdots & O & O\\
O & -s_{1}A_{1} & O & \cdots & O & O\\
O & O & -s_{2}A_{2} & \cdots & O & O\\
\vdots & \vdots & \vdots &  & \vdots & \vdots\\
O & O & O & \cdots & O & -s_{\rho}A_{\rho}
\end{bmatrix}.\label{eq:gamma0'}
\end{equation}
To deal with the rank of $\bar{\Gamma}_{0}$, suppose $(\tau_{n},\tau_{n+1})\in\mathcal{E}^{nb}$,
first one could check that 
\begin{equation}
{\bf span}(\bar{D}({\bf 1}_{\rho}\otimes B_{\hat{A}}))\subset{\bf null}(\Gamma_{0})\label{eq:gamma0-null}
\end{equation}
where ${\bf null}(\hat{A})={\bf span}(\mathcal{B}_{\mathcal{K}_{l}})\cap{\bf span}(\mathcal{B}_{\mathcal{K}_{m}})\cap{\bf null}(A_{n})$.
$\bar{D}$ denotes the NBS's partition of $\{\tau_{1},\tau_{2},...,\tau_{\rho}\}$
which is not respected by the signs of their connecting weights at
$A_{n}$. This demands $A_{n}(v+v)={\bf 0}$ for $\bar{D}({\bf 1}_{\rho}\otimes v)\in{\bf null}(\Gamma_{0})$,
thus the convergence space for the agents will further take intersection
with ${\bf null}(A_{n})$. Here we have
\[
rank(\Gamma_{0})\leqslant(\rho+1)d-dim({\bf null}(\hat{A})).
\]
Instead of directly computing $rank(\bar{\Gamma}_{0})$ with \eqref{eq:rank_block},
we now take a different approach by identifying 
\begin{equation}
{\bf null}\left(\left[\begin{array}{c}
X\\
Y
\end{array}\right]\right)={\bf null}(X)\cap{\bf null}(Y)\label{eq:null_block}
\end{equation}
for block matrix $\begin{bmatrix}X^{T} & Y^{T}\end{bmatrix}^{T}$,
and partitioning $\bar{\Gamma}_{0}$ as $\bar{\Gamma}_{0}=\begin{bmatrix}R'^{T} & Q'^{T}\end{bmatrix}^{T}$
where $R'=\left[\begin{array}{ccccc}
2I & \alpha_{1}I & \cdots & \alpha_{\rho-1}I & -\alpha_{\rho}I\end{array}\right]$, $Q'={\bf blkdiag}\{\bar{A},-s_{1}A_{1},...,-s_{\rho}A_{\rho}\}$.
The eigenvectors that span ${\bf null}(Q')$ are readily obtained
as the columns of 
\begin{align}
\begin{bmatrix}B_{\bar{A}}^{T} & O & \cdots & O\end{bmatrix}^{T},\nonumber \\
\begin{bmatrix}O & O & \cdots & \underset{i\text{-th}}{\underbrace{B_{A_{i}}^{T}}} & \cdots & O\end{bmatrix}^{T}, & i\in\underline{\rho}\label{eq:gamma0'-base}
\end{align}
within which we can find eigenvectors for ${\bf null}(R')$ that are
$\left[\begin{array}{cccccccc}
-v^{T} & {\bf 0} & \cdots & {\bf 0} & 2\alpha_{n}v^{T} & {\bf 0} & \cdots & {\bf 0}\end{array}\right]^{T}$ for $v\in{\bf null}(\hat{A})$. As of now, we have shown the dimension
of ${\bf null}(\Gamma_{0})$ to be no less than $dim({\bf null}(\hat{A}))$
as in eqn. \eqref{eq:gamma0-null}. Suppose ${\bf null}(R')$ has
any other eigenvector as a linear combination of \eqref{eq:gamma0'-base},
i.e., $\begin{bmatrix}k_{0}v_{0}^{T} & k_{1}v_{1}^{T} & \cdots & k_{t}v_{t}^{T} & \cdots & k_{\rho}v_{\rho}^{T}\end{bmatrix}^{T}\in{\bf null}(R')$
where $v_{0}\in{\bf span}(\mathcal{B}_{\bar{A}})\backslash{\bf span}(\mathcal{B}_{\hat{A}})$,
$v_{t}\in{\bf span}(\mathcal{B}_{\bar{A}})\backslash{\bf span}(\mathcal{B}_{\hat{A}})$,
$t\in\underline{\rho}\backslash\{n\}$, $v_{i}\in{\bf span}(\mathcal{B}_{A_{i}}),$
$i\in\underline{\rho}$, and $k_{0},k_{1},...,k_{\rho}$ are not all
zeros. Then Conditions (3) and (4) guarantee that such vectors do
not exist. We have now proven that eqn. \eqref{eq:gamma0-null} is
in fact 
\[
{\bf span}(\bar{D}({\bf 1}_{\rho}\otimes B_{\hat{A}}))={\bf null}(\Gamma_{0}).
\]

\noindent \uuline{Scenario 2}:

We now consider multiple semidefinite paths between $\mathcal{K}_{l}$
and $\mathcal{K}_{m}$, i.e. $\mathring{\mathcal{P}}_{r}=\{(\tau_{1}^{r},\tau_{2}^{r}),...,(\tau_{\rho_{r}}^{r},\tau_{\rho_{r}+1}^{r})\},r\in\underline{\mu}$,
eqn. \eqref{eq:gamma0} is evaluated for every $\mathring{\mathcal{P}}_{r}$
while the agents of the continents $\tau_{1}^{r},\tau_{\rho_{r}+1}^{r}$
are further bound by the bipartite relation $x_{\tau_{1}^{r}}=\pm x_{\tau_{1}^{r'}}$,
$x_{\tau_{1}^{r}}=\pm x_{\tau_{\rho_{r}+1}^{r}}$, $\{r,r'\}\subseteq\underline{\mu}$.
It is then convenient to formulate these relations with the block
matrix 
\begin{equation}
\begin{bmatrix}\Gamma_{1}^{T} & \Gamma_{2}^{T} & \cdots & \Gamma_{\mu}^{T} & \Gamma_{\mu+1}^{T}\end{bmatrix}^{T}\bar{x}={\bf 0}\label{eq:multi-path}
\end{equation}
for all the agents of $\mathcal{K}_{m},\mathcal{K}_{l},$ and $\mathring{\mathcal{P}}_{r},r\in\underline{\mu}$,
where $\Gamma_{1},...,\Gamma_{\mu}$ correspond to the constraints
imposed by $\mathring{\mathcal{P}}_{r},r\in\underline{\mu}$, and
$\Gamma_{\mu+1}$ corresponds to the bipartite relation that couples
the equations of the semidefinite paths. The overall solution of eqn.
\eqref{eq:multi-path} is then found in the intersection of the null
spaces of the blocks. As we have derived ${\bf null}(\Gamma_{0})={\bf span}(\bar{D}({\bf 1}_{\rho}\otimes\bar{B}))$
where $\bar{B}=B_{\bar{A}}$ or $\bar{B}=B_{\hat{A}}$ for each semidefinite
path, and there is Condition (1) to guarantee the solution of the
continents, it is easy to show that the agents of the continents $\mathcal{K}_{l,m}$
and the semidefinite paths must converge bipartitely.

For the entire system $\mathcal{G}$ where several continents are
interconnected through the semidefinite paths, consider each pair
of continents $(\mathcal{K}_{l},\mathcal{K}_{m})$ that sets up a
constraint $\Gamma_{k}^{*}x={\bf 0}$ which eventually stacks up as
$\begin{bmatrix}\Gamma_{1}^{*T} & \Gamma_{2}^{*T} & \cdots & \Gamma_{\theta}^{*T}\end{bmatrix}^{T}x={\bf 0}$.
Similarly, we end up taking the intersection of the solution space
for the set of constraints, and bipartite consensus is achieved due
to the connectivity of the network. In this process, the intersection
is taken with ${\bf null}(A_{e})$ for all $e\in\mathcal{E}^{nb}(\mathcal{V}_{1},\mathcal{V}_{2})$
which means the agents' states converge to ${\bf null}(\mathcal{E}^{nb})$. 
\end{IEEEproof}
\begin{claim}
\label{claim:1}Suppose Assumptions 1 holds. If Conditions (1) and
(2) hold, Conditions (3) and (4) must be satisfied to achieve bipartite
consensus. 
\end{claim}
\begin{IEEEproof}
This can be validated by finding solutions other than the bipartite
one in the Laplacian null space when Condition (3) or (4) does not
hold, e.g., assume ${\bf 0}\neq v\in{\bf null}(A_{i})\cap{\bf null}(A_{j})$,
then $\begin{bmatrix}R^{T} & Q^{T}\end{bmatrix}^{T}$ of eqn. \eqref{eq:gamma0_block}
does not have full rank as $\begin{bmatrix}{\bf 0} & \cdots & \alpha_{i}v^{T} & {\bf 0} & \cdots & -\alpha_{j}v^{T} & {\bf 0} & \cdots & {\bf 0}\end{bmatrix}^{T}\in{\bf null}\left(\begin{bmatrix}R^{T} & Q^{T}\end{bmatrix}^{T}\right)$;
applying the row operations inverse to what transforms eqn. \eqref{eq:gamma0}
to eqn. \eqref{eq:gamma0_block} on this solution, it does not turn
out to be bipartite since there are row blocks of both ${\bf 0}$
and $v$. Notice that the agents in the middle of the semidefinite
paths are not shared with each other, thus even after taking intersection
of the solution spaces of the semidefinite paths, this non-bipartite
solution does not disappear and stays in the Laplacian null space.
The same reason applies to $\bar{\Gamma}_{0}$ of eqn. \eqref{eq:gamma0'}
when Condition (3) does not hold, considering 
\[
\begin{bmatrix}{\bf 0} & {\bf 0} & \cdots & \alpha_{i}v^{T} & \cdots & -\alpha_{j}v^{T} & \cdots & {\bf 0}\end{bmatrix}^{T}\in{\bf null}(\bar{\Gamma}_{0}).
\]
If Condition (4) does not hold, there exist $v_{0}\in{\bf span}(B_{\bar{A}})\backslash{\bf span}(B_{\hat{A}})$,
$v_{i}\in{\bf span}(B_{A_{i}})$, $i\in\underline{\rho}\backslash\{n\}$,
such that
\[
\begin{bmatrix}v_{0}^{T} & v_{1}^{T} & \cdots & v_{n-1}^{T} & {\bf 0} & v_{n+1}^{T} & \cdots & v_{\rho}^{T}\end{bmatrix}^{T}\in{\bf null}(\bar{\Gamma}_{0}),
\]
due to which bipartite consensus is not achieved for agents on the
semidefinite path.
\end{IEEEproof}
\bibliographystyle{IEEEtran}
\bibliography{2_Users_chloekiss_Dropbox_Matrix-bipartite-cluster_Tech_note_suff_clus}

\begin{thebibliography}{10}
\providecommand{\url}[1]{#1}
\csname url@samestyle\endcsname
\providecommand{\newblock}{\relax}
\providecommand{\bibinfo}[2]{#2}
\providecommand{\BIBentrySTDinterwordspacing}{\spaceskip=0pt\relax}
\providecommand{\BIBentryALTinterwordstretchfactor}{4}
\providecommand{\BIBentryALTinterwordspacing}{\spaceskip=\fontdimen2\font plus
\BIBentryALTinterwordstretchfactor\fontdimen3\font minus
  \fontdimen4\font\relax}
\providecommand{\BIBforeignlanguage}[2]{{%
\expandafter\ifx\csname l@#1\endcsname\relax
\typeout{** WARNING: IEEEtran.bst: No hyphenation pattern has been}%
\typeout{** loaded for the language `#1'. Using the pattern for}%
\typeout{** the default language instead.}%
\else
\language=\csname l@#1\endcsname
\fi
#2}}
\providecommand{\BIBdecl}{\relax}
\BIBdecl

\bibitem{Castellano2009}
C.~Castellano, S.~Fortunato, and V.~Loreto, ``Statistical physics of social
  dynamics,'' \emph{Rev. Mod. Phys.}, vol.~81, pp. 591--646, May 2009.

\bibitem{Proskurnikov2016}
A.~V. Proskurnikov, A.~S. Matveev, and M.~Cao, ``Opinion dynamics in social
  networks with hostile camps: Consensus vs. polarization,'' \emph{IEEE
  Transactions on Automatic Control}, vol.~61, no.~6, pp. 1524--1536, 2016.

\bibitem{Doerfler2013}
F.~D\"{o}rfler, M.~Chertkov, and F.~Bullo, ``Synchronization in complex
  oscillator networks and smart grids,'' \emph{Proceedings of the National
  Academy of Sciences}, vol. 110, no.~6, pp. 2005--2010, 2013.

\bibitem{Kamel2020}
M.~A. Kamel, X.~Yu, and Y.~Zhang, ``Formation control and coordination of
  multiple unmanned ground vehicles in normal and faulty situations: A
  review,'' \emph{Annual Reviews in Control}, vol.~49, pp. 128--144, 2020.

\bibitem{Tahir2019}
A.~Tahir, J.~B{\"o}ling, M.-H. Haghbayan, H.~T. Toivonen, and J.~Plosila,
  ``Swarms of unmanned aerial vehicles a survey,'' \emph{Journal of Industrial
  Information Integration}, vol.~16, p. 100106, 2019.

\bibitem{Liu2023}
G.~Liu, L.~Chen, K.~Liu, and Y.~Luo, ``A swarm of unmanned vehicles in the
  shallow ocean: A survey,'' \emph{Neurocomputing}, vol. 531, pp. 74--86, 2023.

\bibitem{olfati2004consensus}
R.~Olfati-Saber and R.~M. Murray, ``Consensus problems in networks of agents
  with switching topology and time-delays,'' \emph{IEEE Transactions on
  Automatic Control}, vol.~49, no.~9, pp. 1520--1533, 2004.

\bibitem{altafini2012consensus}
C.~Altafini, ``Consensus problems on networks with antagonistic interactions,''
  \emph{IEEE Transactions on Automatic Control}, vol.~58, no.~4, pp. 935--946,
  2012.

\bibitem{jadbabaie2003coordination}
A.~Jadbabaie, J.~Lin, and A.~S. Morse, ``Coordination of groups of mobile
  autonomous agents using nearest neighbor rules,'' \emph{IEEE Transactions on
  Automatic Control}, vol.~48, no.~6, pp. 988--1001, 2003.

\bibitem{tuna2016synchronization}
S.~E. Tuna, ``Synchronization under matrix-weighted laplacian,''
  \emph{Automatica}, vol.~73, pp. 76--81, 2016.

\bibitem{TRINH2018415}
M.~H. Trinh, C.~Van Nguyen, Y.-H. Lim, and H.-S. Ahn, ``Matrix-weighted
  consensus and its applications,'' \emph{Automatica}, vol.~89, pp. 415 -- 419,
  2018.

\bibitem{Pan2019}
L.~Pan, H.~Shao, M.~Mesbahi, Y.~Xi, and D.~Li, ``Bipartite consensus on
  matrix-valued weighted networks,'' \emph{IEEE Transactions on Circuits and
  Systems II: Express Briefs}, vol.~66, no.~8, pp. 1441--1445, 2019.

\bibitem{Ramirez2009}
J.~L. Ramirez, M.~Pavone, E.~Frazzoli, and D.~W. Miller, ``Distributed control
  of spacecraft formation via cyclic pursuit: Theory and experiments,'' in
  \emph{2009 American Control Conference}, 2009, pp. 4811--4817.

\bibitem{zhao2016localizability}
S.~Zhao and D.~Zelazo, ``Localizability and distributed protocols for
  bearing-based network localization in arbitrary dimensions,''
  \emph{Automatica}, vol.~69, pp. 334--341, 2016.

\bibitem{Lee2016}
B.-H. Lee and H.-S. Ahn, ``Distributed formation control via global orientation
  estimation,'' \emph{Automatica}, vol.~73, pp. 125--129, 2016.

\bibitem{Hong2011}
H.~Hong and S.~H. Strogatz, ``Kuramoto model of coupled oscillators with
  positive and negative coupling parameters: An example of conformist and
  contrarian oscillators,'' \emph{Phys. Rev. Lett.}, vol. 106, p. 054102, Feb
  2011.

\bibitem{Hu2015}
J.~Hu and H.~Zhu, ``Adaptive bipartite consensus on coopetition networks,''
  \emph{Physica D: Nonlinear Phenomena}, vol. 307, pp. 14--21, 2015.

\bibitem{Zong2019}
C.~Zong, Z.~Ji, L.~Tian, and Y.~Zhang, ``Distributed multi-robot formation
  control based on bipartite consensus with time-varying delays,'' \emph{IEEE
  Access}, vol.~7, pp. 144\,790--144\,798, 2019.

\bibitem{su2019bipartite}
H.~Su, J.~Chen, Y.~Yang, and Z.~Rong, ``The bipartite consensus for multi-agent
  systems with matrix-weight-based signed network,'' \emph{IEEE Transactions on
  Circuits and Systems II: Express Briefs}, 2019.

\bibitem{wang2022characterizing}
C.~Wang, L.~Pan, H.~Shao, D.~Li, and Y.~Xi, ``Characterizing bipartite
  consensus on signed matrix-weighted networks via balancing set,''
  \emph{Automatica}, vol. 141, p. 110237, 2022.

\bibitem{tuna2019synchronization}
S.~E. Tuna, ``Synchronization of small oscillations,'' \emph{Automatica}, vol.
  107, pp. 154--161, 2019.

\bibitem{Ahn2020}
H.-S. Ahn, Q.~V. Tran, M.~H. Trinh, M.~Ye, J.~Liu, and K.~L. Moore, ``Opinion
  dynamics with cross-coupling topics: Modeling and analysis,'' \emph{IEEE
  Transactions on Computational Social Systems}, vol.~7, no.~3, pp. 632--647,
  2020.

\end{thebibliography}

\end{document}